\documentclass[a4paper,11pt,reqno]{amsart}
\usepackage[T1]{fontenc}
\usepackage{lmodern}
\usepackage[a4paper,margin=27mm]{geometry}
\usepackage{amsmath,amssymb,amsthm,mathtools}
\usepackage{mathrsfs}
\usepackage{xcolor}
\usepackage{tikz}
\usetikzlibrary{shapes.geometric}
\usepackage{placeins}
\usepackage{enumitem,needspace}
\usepackage{hyperref}
\usepackage{xurl}
\hypersetup{colorlinks=true,linkcolor=black,citecolor=black,
 urlcolor=black,pdftitle={Polypositroids are Ehrhart positive},
 pdfauthor={Yao Li and Sihong Shao}}
\numberwithin{equation}{section}
\newtheorem{theorem}{Theorem}[section]
\newtheorem{proposition}[theorem]{Proposition}
\newtheorem{lemma}[theorem]{Lemma}

\theoremstyle{definition}
\newtheorem{definition}[theorem]{Definition}
\newtheorem{question}[theorem]{Question}
\newtheorem{example}{Example}[section]
\theoremstyle{remark}
\newtheorem{remark}[theorem]{Remark}
\newcommand{\R}{\mathbb R}
\newcommand{\Z}{\mathbb Z}

\allowdisplaybreaks[2]
\setlist{itemsep=0.2em,topsep=0.4em}
\title[Polypositroids are Ehrhart positive]{Polypositroids are Ehrhart positive}
\author{YAO LI AND SIHONG SHAO}
\date{}
\begin{document}
\begin{abstract}
We prove that integral polypositroids are Ehrhart positive.
In particular, all positroid base
polytopes are Ehrhart positive. An inclusion-exclusion recurrence for lattice-point
generating functions is derived from the interval inequalities.
The generating functions are represented by operators whose integral kernels
are computed using divided differences. Weight-preserving cancellation reduces
the proof to two endpoint corrections given by one-sided kernel limits
at the lower and upper bounds on the coordinate sum. A volume estimate
near each endpoint shows that both corrections have nonnegative coefficients.
The Ehrhart polynomial is then recovered from the lower correction at the
prescribed total coordinate sum.
\end{abstract}
\maketitle

\section{Introduction}\label{sec:introduction}

Let $E$ be a finite set. For an \emph{integral} polytope
$\mathcal P\subseteq\R^E$, that is, a polytope with vertices in
$\Z^E$, define
\[
 i(\mathcal P,t)\coloneqq\#\bigl(t\mathcal P\cap\Z^E\bigr).
\]
Unless otherwise stated, the dilation parameter $t$ is a nonnegative integer.
By Ehrhart's theorem \cite[Theorem~3.8]{BeckRobins},
$i(\mathcal P,t)$ is a polynomial in $t$ of degree
$\dim\mathcal P$ with rational coefficients, called the
\emph{Ehrhart polynomial} of $\mathcal P$. We say that $\mathcal P$
is \emph{Ehrhart positive} if
all coefficients of $i(\mathcal P,t)$ are nonnegative.

De Loera, Haws, and K\"oppe conjectured that matroid base polytopes
are Ehrhart positive
\cite[Conjecture~2(B)]{DHK}.
For $i\in E$, let $\mathbf e_i$ denote the $i$th standard unit
vector of $\R^E$. A \emph{generalized permutohedron} is a nonempty
polytope in $\R^E$ whose edges are parallel to vectors
$\mathbf e_i-\mathbf e_j$ for distinct $i,j\in E$.
Equivalently, generalized permutohedra are the base polytopes of
submodular functions \cite[Theorem~4.1]{LP}.
Matroid base polytopes form a subclass of generalized permutohedra.
Castillo and Liu proposed the broader conjecture that every integral
generalized permutohedron is Ehrhart positive
\cite[Conjecture~1.2]{CL}.
Among generalized permutohedra, Ehrhart positivity is known for
nonnegative integral Minkowski sums of coordinate simplices
$\operatorname{conv}\{\mathbf e_i:i\in A\}$, called
$\mathcal Y$-generalized permutohedra
\cite[Theorem~11.3]{PostnikovGP}, and for integral generalized
permutohedra of dimension at most six \cite{CL}.
Ferroni nevertheless constructed matroid base polytopes with negative
Ehrhart coefficients, disproving both conjectures \cite{Ferroni}.
The Ehrhart-positivity conjecture for Schubitopes proposed by
Monical, Tokcan, and Yong \cite{MTY}
was also disproved by Li and St.~Dizier \cite{LiStDizier}.

\emph{Positroids} are matroids represented by real matrices of full
row rank with nonnegative maximal minors. Positroids index the
cells of the totally nonnegative Grassmannian \cite{Postnikov}.
The study of positroids connects Grassmannian geometry and topology
\cite{GKL,KLS}, cluster algebras \cite{GL}, and mathematical physics
\cite{AHBCGPT}.
Ferroni, Jochemko, and Schr\"oter conjectured that positroid base
polytopes are Ehrhart positive \cite[Conjecture~6.3]{FJS}.

Ehrhart positivity has been proved for several classes of matroids
that become positroids after relabeling.
Early results include Ferroni's proofs for uniform
\cite{FerroniUniform} and minimal matroids \cite{FerroniMinimal},
and the rank-two result of Ferroni, Jochemko, and Schr\"oter \cite{FJS}.
Further progress came from Fan and Li's work on certain families
of Schubert matroids \cite{FanLi} and from Deligeorgaki, McGinnis,
and Vindas-Mel\'endez's proof for panhandle matroids \cite{DMV}.
Ferroni, Morales, and Panova proved positivity for snake matroids
\cite{FMPSkew}, while Chen, Li, and Yao established positivity for
$(a,b)$-Catalan matroids \cite{ChenLiYao}.
Ferroni, Morales, and Panova later proved positivity for all
lattice path matroids \cite{FMP}.
More recently, Li and Shao proved positivity for series-parallel
extensions of multi-path matroids \cite{LiShaoSP}, and Liu and Zhang
for paving positroids \cite{LiuZhang}.

Equip $E$ with a cyclic order. A \emph{cyclic interval} is a nonempty
subset consisting of consecutive elements in this cyclic order;
we allow $E$ itself as a cyclic interval.
For a coordinate vector $x=(x_i)_{i\in E}$ and $A\subseteq E$,
write $x(A)\coloneqq\sum_{i\in A}x_i$. We use the same sum notation
for other coordinate vectors.
A \emph{type-$A$ alcoved polytope} is a polytope in $\R^E$
defined by an equation $x(E)=r$ and inequalities $x(I)\le a_I$
for cyclic intervals $I\subseteq E$, where $r$ and each $a_I$
are real constants \cite[Section~3.2]{LP}.
The class of \emph{polypositroids}, introduced by Lam and Postnikov,
is exactly the intersection of the classes of generalized
permutohedra and type-$A$ alcoved polytopes
\cite[Definition~3.8]{LP}.
Unless stated otherwise, the polypositroids considered below are integral.
Positroid base polytopes are polypositroids
\cite[Theorem~2.1]{LP}.
The theory of polypositroids also has connections with Coxeter
combinatorics and cluster algebras
\cite[Sections~1.1 and~1.3]{LP}.

Earlier positive cases include integral Pitman--Stanley polytopes,
viewed as polypositroids up to integral affine equivalence
\cite[Corollary~10]{StanleyPitman}.
Ferroni and McGinnis later proved positivity for another subclass
of polypositroids: intersections of boxes with integer coordinate
bounds and hyperplanes $x(E)=r$, where $r$ is an integer \cite{FMPrisms}.
The Ehrhart-positivity conjecture for polypositroids was
explicitly suggested by Jiang \cite[Section~5.7, p.~66]{Jiang}.
We prove the following theorem.

\begin{theorem}\label{thm:main}
Integral polypositroids are Ehrhart positive.
\end{theorem}

In particular, Theorem~\ref{thm:main} settles the conjecture of
Ferroni, Jochemko, and Schr\"oter \cite[Conjecture~6.3]{FJS}
that positroids are Ehrhart positive.

Postnikov introduced \emph{flowtroid polytopes} using directed
networks with designated source vertices
\cite[pp.~543--545]{PostnikovFlowtroids}.
Allow nonnegative flows on the edges. At each vertex outside the
designated source set, require the total inflow minus the total
outflow to equal a prescribed integer. The net outflow vectors
at the designated sources realized by edge flows satisfying these
constraints form a flowtroid polytope.
Flowtroid polytopes are integral generalized permutohedra and include
both integral polypositroids and $\mathcal Y$-generalized permutohedra
\cite[p.~544]{PostnikovFlowtroids}.
In comparing these classes, we allow coordinate permutations and
integer translations, both of which preserve Ehrhart polynomials.

Among matroid base polytopes, the flowtroid realizations are exactly
those of gammoids \cite[p.~545]{PostnikovFlowtroids}.
A gammoid is a matroid on a subset of the vertices of a directed graph;
a set is independent exactly when its elements can be linked to a
fixed target set by vertex-disjoint directed paths.
Ferroni's sparse paving counterexample \cite[Theorem~5.3]{Ferroni}
has an $M(K_4)$-minor and hence is not a gammoid, so its base
polytope lies outside the flowtroid class.\footnote{Let $M$ be the
matroid from Ferroni's $C_0$ construction, with ground set labeled
by $0,\ldots,19$;
the nine-element nonbases are precisely the subsets whose labels
sum to $0$ modulo $20$. Contract $C=\{2,4,5,8,9,18\}$ and
restrict to $D=\{0,1,3,10,11,13\}$. The correspondence
$(0,1,3,13,11,10)\mapsto(12,13,14,23,24,34)$ with the edges
of $K_4$ identifies $(M/C)|D$ with the cycle matroid $M(K_4)$.
Gammoids are minor-closed and base-orderable, whereas $M(K_4)$
is not base-orderable \cite[Section~1]{BoninSavitsky}.}
Every transversal matroid is a gammoid \cite{IngletonPiff}.
The connected transversal matroid of rank~$45$ on $61$ elements
constructed in \cite[Proposition~7.1]{LiShaoSP} has a base polytope
with a negative quadratic Ehrhart coefficient.
Consequently, flowtroid polytopes are not all Ehrhart positive.

\begin{figure}[htbp]
\centering
\begin{tikzpicture}[x=1cm,y=1cm,font=\small,
 every node/.style={align=center,inner sep=1pt},
 boundary/.style={draw=black!65,line width=.45pt}]
 \fill[black!7] (3,1) rectangle (8.4,5);
 \fill[black!7] (10.1,2.95) ellipse (3.7 and 1.45);
 \draw[line width=.65pt] (3,0) rectangle (14.7,7.3);
 \draw[boundary] (3,.45) rectangle (14.2,6.05);
 \draw[line width=.65pt] (0,1) rectangle (8.4,5);
 \draw[boundary] (10.1,2.95) ellipse (3.7 and 1.45);
 \node[font=\normalsize] at (8.85,6.97) {Generalized permutohedra};
 \fill (4,6.45) circle (1.8pt);
 \node[anchor=west,font=\footnotesize] at (4.25,6.45)
   {A matroid counterexample \cite{Ferroni}};
 \node[font=\normalsize] at (8.6,5.75) {Flowtroid polytopes};
 \node[text width=2.55cm] at (1.5,3.4)
   {Type-$A$ alcoved\\polytopes};
 \fill (1.5,2.5) circle (1.8pt);
 \node[font=\scriptsize] at (1.5,1.95)
   {\shortstack{Liu--Tsuchiya\\counterexample \cite{LiuTsuchiya}}};
 \node at (5.1,4.55) {Polypositroids};
 \node[font=\scriptsize] at (5.2,4.15)
   {Ehrhart positive: Theorem~\ref{thm:main}};
 \draw[boundary] (5,2.8) ellipse (1.9 and .55);
 \node at (4.9,2.8) {Positroids};
 \draw[boundary] (7.4,2.85) circle (.6);
 \node[font=\fontsize{7}{8}\selectfont] at (7.44,2.85)
   {\shortstack{Pitman--\\Stanley}};
 \node[text width=3.7cm] at (10.7,3.4)
   {$\mathcal Y$-generalized\\permutohedra};
 \node[font=\scriptsize] at (10.7,2.65)
   {Ehrhart positive \cite{PostnikovGP}};
 \fill (13.5,4.85) circle (1.8pt);
 \node[anchor=east,font=\footnotesize,text width=3.6cm] at (13.25,4.85)
   {Transversal matroid\\counterexample \cite{LiShaoSP}};
\end{tikzpicture}
\caption{All polytopes are integral.
Shading indicates Ehrhart-positive classes.}
\label{fig:positivity-classes}
\end{figure}
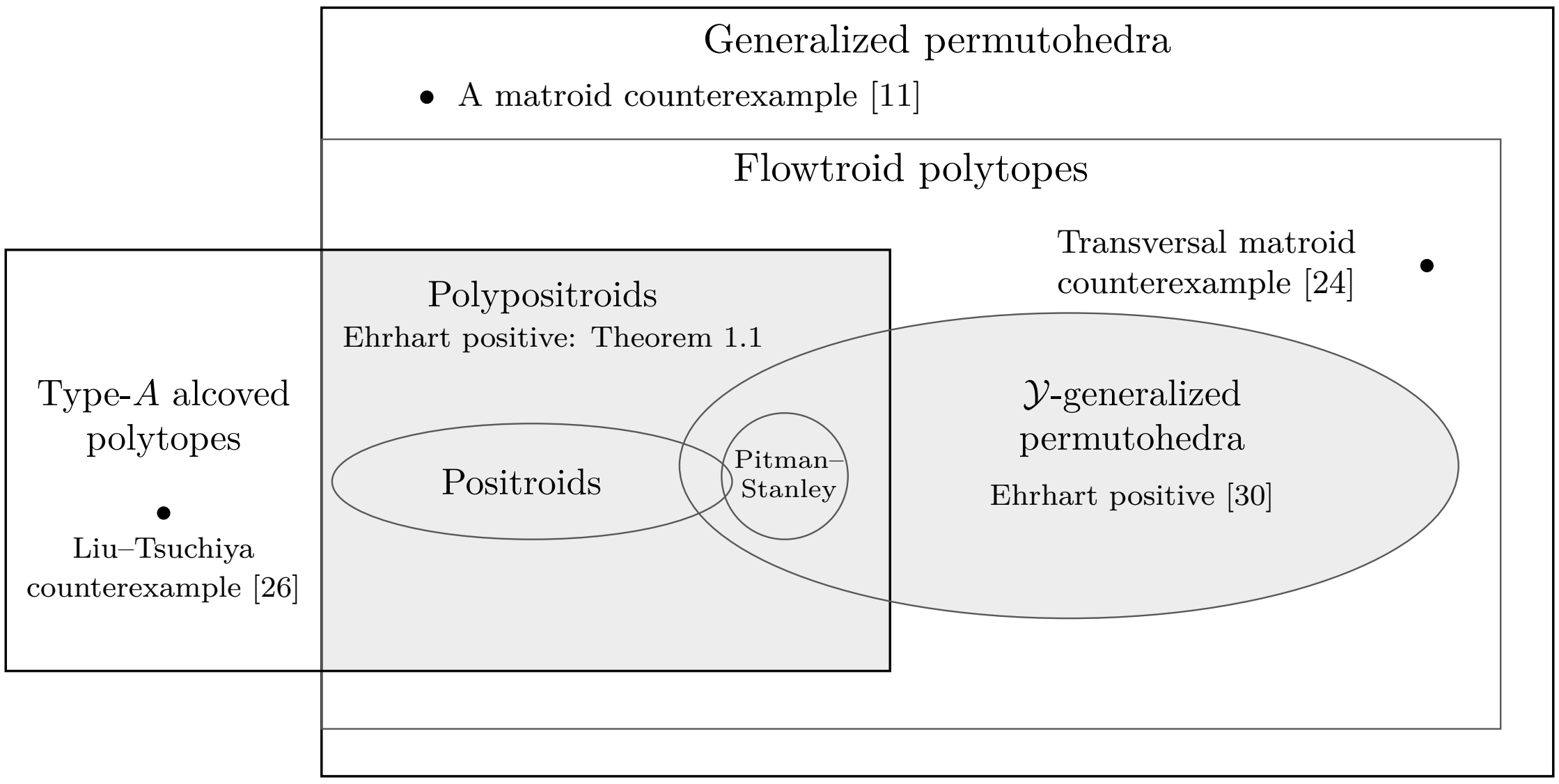

Within the inclusion chain from positroid base polytopes through
integral polypositroids and flowtroid polytopes to integral generalized
permutohedra, integral polypositroids form the largest known
Ehrhart-positive class.
The Ehrhart-positive class $\mathcal Y$ extends beyond polypositroids,
as illustrated in Figure~\ref{fig:positivity-classes}.

Neither integral generalized permutohedra nor integral type-$A$
alcoved polytopes are Ehrhart positive in general
\cite{Ferroni,LiuTsuchiya}. Theorem~\ref{thm:main} establishes
Ehrhart positivity for their intersection, the class of integral
polypositroids. This leads to the following fundamental question.

\begin{question}\label{ques:positive-classification}
What geometric or combinatorial conditions are necessary and
sufficient for an integral polytope to be Ehrhart positive?
\end{question}

Within generalized permutohedra, a step toward
Question~\ref{ques:positive-classification} would be to find a
common structural sufficient condition encompassing polypositroids
and $\mathcal Y$-generalized permutohedra.

\FloatBarrier

Our proof combines interval inclusion-exclusion with coefficientwise
estimates for integral operators. Section~\ref{sec:geometry} first
reduces the problem to a polypositroid of dimension $|E|-1$ whose
coordinate minima are zero. After deleting one coordinate, we work
with projections onto intervals $I$ in the remaining linear
order. For a positive integer $t$, the polynomial $\beta_I(t,z)$ records
each lattice point $a$ in the $t$-fold dilate of the projection
by the monomial $z^{a(I)}$.
Writing $\underline b(I)$ and $\overline b(I)$ for the minimum and maximum of $x(I)$
on $\mathcal P$, the allowed total coordinate sums satisfy
\[
 t\underline b(I)\le a(I)\le t\overline b(I).
\]
The endpoint equalities $a(I)=t\underline b(I)$ and $a(I)=t\overline b(I)$ correspond
to terms with exponents $t\underline b(I)$ and $t\overline b(I)$ in the generating
function $\beta_I(t,z)$ defined in
\eqref{eq:interval-generating-function}. Both endpoint terms are retained.

Section~\ref{sec:counting} derives an inclusion-exclusion
recurrence indexed by partitions $\pi$ of $I$ into at least two
intervals. A partition contributes
$(-1)^{|\pi|}\prod_{J\in\pi}\beta_J(t,z)$, before retaining
only exponents in $[t\underline b(I),\,t\overline b(I)]$.
The recurrence counts lattice points for each positive integer $t$,
but does not by itself determine the signs of the Ehrhart coefficients.
To study those coefficients, Section~\ref{sec:operators} represents
the generating functions by operators that combine values of a test function at specified
points with integrals of that function. Divided differences express
the integral kernels as polynomials in $t$.
The lower and upper coordinate-sum bounds become restrictions on
the accumulated change in the test-function argument; the separate
evaluation terms at both endpoints remain part of the representation.

The structural assertion stated in Section~\ref{sec:recursion}
and proved in Section~\ref{sec:main-proof} is that the operator
for every interval $I$ has evaluation weights
and an integral kernel with nonnegative coefficients in $t$.
The proof pairs terms from different interval partitions by adding or
removing a cut between adjacent blocks. The paired terms have the same
polynomial weight and integration measure and opposite signs.
For fixed blocks, operator branches, and displacements, the signed
multiplicity is zero or one whenever the displacements sum to a value in
$[\underline b(I),\overline b(I)]$. At distance
$s$ just outside an endpoint, the submodular interval inequalities bound
the total measure of the regions of negative multiplicity by $O(s)$.
In each degree in $t$, the integrands are uniformly bounded, so their
integrals over these regions tend to zero as $s\downarrow0$.
The kernel limits give the corrections added to the original evaluation
weights at the lower and upper endpoints. Section~\ref{sec:main-proof}
proves nonnegativity of the corrections by induction on $|I|$ and uses
the source expansion to establish positivity, support, and the endpoint
conditions for the interval operators.
The Ehrhart polynomial is then identified with a coefficientwise
limit of nonnegative kernel polynomials, proving
Theorem~\ref{thm:main}.

\FloatBarrier
\section{Lattice-point generating functions}\label{sec:geometry}

Let $\mathcal P\subseteq\{x(E)=r\}$ be a polypositroid.
For each $A\subseteq E$, define the \emph{tight coordinate-sum bounds}
for $x(A)$ on $\mathcal P$ by
\begin{equation}\label{eq:1-1}
 \overline b(A)\coloneqq\max_{x\in \mathcal P}x(A),\qquad
 \underline b(A)\coloneqq\min_{x\in \mathcal P}x(A)=r-\overline b(E\setminus A).
\end{equation}
Since $\mathcal P$ is integral, $\overline b(A),\underline b(A)\in\Z$ for every $A\subseteq E$;
also $\overline b(\varnothing)=\underline b(\varnothing)=0$.
For a generalized permutohedron, the tight coordinate-sum upper
bounds form a submodular function and give the inequality
description \cite[Theorem~4.1]{LP}
\begin{equation}\label{eq:full-support}
 \mathcal P=\{x\in\R^E:x(E)=r,\ x(A)\le \overline b(A)
                         \text{ for every }A\subseteq E\}.
\end{equation}
Conversely, any submodular function $\overline b:2^E\to\R$ with
$\overline b(\varnothing)=0$ and $\overline b(E)=r$ defines a generalized permutohedron
by \eqref{eq:full-support}, called the \emph{base polytope} of $\overline b$.
If $\overline b$ is integer-valued, the base polytope is integral, since each
vertex coordinate $x_i$ is a difference $\overline b(A\cup\{i\})-\overline b(A)$
for some $A\subseteq E\setminus\{i\}$.
The fixed polypositroid $\mathcal P$ is also alcoved, so $\mathcal P$ can be defined
inside $\{x(E)=r\}$ using only cyclic-interval inequalities
$x(J)\le a_J$. For each defining inequality, the maximum
$\overline b(J)$ satisfies $\overline b(J)\le a_J$. Every point of $\mathcal P$
satisfies $x(J)\le \overline b(J)$, while any point satisfying the tighter
bound also satisfies $x(J)\le \overline b(J)\le a_J$.
Replacing the defining bounds by $\overline b(J)$, and including the tight
bounds for the remaining cyclic intervals, therefore gives
\[
 \mathcal P=\{x\in\R^E:x(E)=r,\ x(J)\le \overline b(J)
                 \text{ for every cyclic interval }J\subseteq E\}.
\]
The function $\overline b$ remains defined on every subset of $E$ by
\eqref{eq:1-1}; the constraints for noncyclic subsets need not be
listed in this description of $\mathcal P$.
Submodularity of $\overline b$ and the identity $\underline b(A)=r-\overline b(E\setminus A)$
imply that $\underline b$ is supermodular.
Submodularity of $\overline b$ also gives, for all $A,B\subseteq E$,
\begin{equation}\label{eq:1-4}
 \overline b(A)-\underline b(B)\ge \overline b(A\setminus B)-\underline b(B\setminus A).
\end{equation}

Let $A,B\subseteq E$ be intersecting intervals in a
fixed linear order on $E$, with neither interval containing the other.
Submodularity of $\overline b$, supermodularity of $\underline b$, and \eqref{eq:1-4} give
\begin{equation}\label{eq:local-crossing}
\begin{aligned}
 \overline b(A)+\overline b(B)&\ge \overline b(A\cup B)+\overline b(A\cap B),\\
 \underline b(A)+\underline b(B)&\le \underline b(A\cup B)+\underline b(A\cap B),\\
 \overline b(A)-\underline b(B)&\ge \overline b(A\setminus B)-\underline b(B\setminus A),\\
 \overline b(B)-\underline b(A)&\ge \overline b(B\setminus A)-\underline b(A\setminus B).
\end{aligned}
\end{equation}
The sets $A\cup B$, $A\cap B$, $A\setminus B$, and $B\setminus A$
are nonempty intervals. Section~\ref{sec:violations}
uses \eqref{eq:local-crossing} to control intersecting intervals
whose coordinate sums violate their bounds.

For $A\subseteq B\subseteq E$, let
\[
 \operatorname{pr}_A:\R^B\longrightarrow\R^A,\qquad
 x\longmapsto(x_i)_{i\in A},
\]
be the \emph{coordinate projection} onto $\R^A$.
For a partition $E=C_1\sqcup\cdots\sqcup C_s$ and each
$a\in\{1,\ldots,s\}$, put
\[
 \mathcal P_a\coloneqq\operatorname{pr}_{C_a}(\mathcal P).
\]
The elements of each $C_a$ retain their relative cyclic order in $E$.

\begin{proposition}\label{prop:components}
There is a partition $E=C_1\sqcup\cdots\sqcup C_s$ such that
\[
 \mathcal P=\prod_{a=1}^s \mathcal P_a.
\]
Each factor $\mathcal P_a$ is a polypositroid on $C_a$ of dimension
$|C_a|-1$, and
\begin{equation}\label{eq:ehrhart-product}
 i(\mathcal P,t)=\prod_{a=1}^s i(\mathcal P_a,t).
\end{equation}
\end{proposition}

\begin{proof}
By \cite[Theorem~6.12, Lemma~7.1, and Proposition~7.2]{LP}, there is
a partition $E=C_1\sqcup\cdots\sqcup C_s$ for which each projection
$\mathcal P_a$ is a polypositroid on $C_a$ and
\[
 \mathcal P=\prod_{a=1}^s\mathcal P_a,\qquad
 \dim\mathcal P=|E|-s.
\]
For each $a$, the factor $\mathcal P_a$ lies in the hyperplane
$x(C_a)=\overline b(C_a)$, so $\dim\mathcal P_a\le |C_a|-1$.
The product decomposition gives
\[
 \sum_{a=1}^s\dim\mathcal P_a
 =\dim\mathcal P
 =|E|-s
 =\sum_{a=1}^s(|C_a|-1).
\]
Therefore, $\dim\mathcal P_a=|C_a|-1$ for every $a$.

Every vertex of a coordinate projection of an integral polytope is
the projection of an integral vertex, so $\mathcal P_a$ is integral.
For every nonnegative integer $t$,
\[
 i(\mathcal P,t)
 =\#\!\left(\prod_{a=1}^s
       (t\mathcal P_a\cap\Z^{C_a})\right)
 =\prod_{a=1}^s\#(t\mathcal P_a\cap\Z^{C_a})
 =\prod_{a=1}^s i(\mathcal P_a,t).
\]
\end{proof}

When $\mathcal P$ is a matroid base polytope, the factorization in
Proposition~\ref{prop:components} is the product decomposition induced
by the connected components of the matroid.

By \eqref{eq:ehrhart-product},
it suffices to prove coefficient nonnegativity for each factor
$\mathcal P_a$. If $|C_a|=1$, then $\dim\mathcal P_a=0$, so
$\mathcal P_a$ consists of a single lattice point and has Ehrhart
polynomial $1$.
Fix a factor $\mathcal P_a$ with $|C_a|\ge2$, and denote it by
$\mathcal P\subseteq\R^E$, where $E=C_a$.
Let $r$ be the constant value of $x(E)$ on $\mathcal P$, and
define $\overline b$ and $\underline b$ for $\mathcal P$ by \eqref{eq:1-1}.
Put $n=|E|$; then $n\ge2$ and $\dim\mathcal P=n-1$.

To make the minimum of each coordinate equal to zero, set
$\ell_i\coloneqq \underline b(\{i\})$ for $i\in E$ and consider
$\mathcal P-\ell$, where $\ell=(\ell_i)_{i\in E}\in\Z^E$.
For $A\subseteq E$, write $\ell(A)\coloneqq\sum_{i\in A}\ell_i$.
Under the map $x\mapsto x-\ell$, the tight bounds $\overline b(A)$ and $\underline b(A)$
become $\overline b(A)-\ell(A)$ and $\underline b(A)-\ell(A)$, respectively, and the
total coordinate sum becomes $r-\ell(E)$.
The polytope $\mathcal P-\ell$ has the same edge directions as
$\mathcal P$ and is defined inside $\{x(E)=r-\ell(E)\}$ by
$x(J)\le \overline b(J)-\ell(J)$ for cyclic intervals $J\subseteq E$.
Thus $\mathcal P-\ell$ is again a polypositroid.
For each nonnegative integer $t$, the map $x\mapsto x-t\ell$
is a bijection from $t\mathcal P\cap\Z^E$ to
$t(\mathcal P-\ell)\cap\Z^E$, so
$i(\mathcal P-\ell,t)=i(\mathcal P,t)$.
For every pair of subsets $A,B\subseteq E$, the coordinate sums satisfy
\[
\begin{aligned}
 \ell(A)+\ell(B)&=\ell(A\cup B)+\ell(A\cap B),\\
 \ell(A)-\ell(B)&=\ell(A\setminus B)-\ell(B\setminus A).
\end{aligned}
\]
Replacing $\overline b,\underline b$ by $\overline b-\ell,\underline b-\ell$ subtracts
$\ell(A)+\ell(B)$ from both sides of the submodularity and
supermodularity inequalities, and $\ell(A)-\ell(B)$ from both
sides of \eqref{eq:1-4}. Hence all four inequalities in
\eqref{eq:local-crossing} remain valid.
Until the final application of Proposition~\ref{prop:components},
write $\mathcal P$ for the translated polypositroid $\mathcal P-\ell$.
Let $r$ be its constant coordinate sum, and define $\overline b,\underline b$ for the
translated polypositroid by \eqref{eq:1-1}.
Translation preserves dimension, so $\dim\mathcal P=n-1$.
We claim that each coordinate function is nonconstant on $\mathcal P$.
Otherwise, $\mathcal P$ would satisfy an additional equation
$x_i=c$, independent of $x(E)=r$ because $n\ge2$.
This would imply $\dim\mathcal P\le n-2$, a contradiction.
The zero coordinate minima and integrality of $\mathcal P$ therefore
give, for every $i\in E$,
\[
 \min_{x\in \mathcal P}x_i=0,\qquad
 c_i\coloneqq\max_{x\in \mathcal P}x_i\in\Z_{\ge1}.
\]
For each $i\in E$, the projection of $\mathcal P$ onto coordinate
$i$ is the interval $[0,c_i]$.
Since every coordinate is nonnegative on $\mathcal P$,
$\underline b(A)\ge0$ for every $A\subseteq E$ and $r\ge c_i\ge1$ for every
$i\in E$.
For any nonempty proper subset $A\subsetneq E$, an equation
$x(A)=c$ would likewise be independent of $x(E)=r$ and force
$\dim\mathcal P\le n-2$. Hence $x(A)$ is nonconstant on
$\mathcal P$, giving
\begin{equation}\label{eq:1-5}
 \overline b(A)-\underline b(A)>0.
\end{equation}
In particular, $\overline b(I)-\underline b(I)$ is a positive integer for every
proper cyclic interval $I\subsetneq E$.

We now use the normalized bounds to describe coordinate projections
and their lattice-point generating functions. Fix $e\in E$ and
put $N\coloneqq E\setminus\{e\}$. Order $N$ by following the
cyclic order on $E$ from the successor of $e$ to its predecessor.
For $x\in\mathcal P$, write $y=\operatorname{pr}_N(x)$. The constraint for
the whole interval $E$ is the equation $x(E)=r$, which determines
the deleted coordinate. Denote the projection by
\begin{equation}\label{eq:1-3}
 \mathcal P_N\coloneqq\operatorname{pr}_N(\mathcal P),\qquad x_e=r-y(N).
\end{equation}
Thus $\operatorname{pr}_N:\mathcal P\to\mathcal P_N$ is a bijection, with
inverse obtained by adjoining $x_e=r-y(N)$.

To describe $\mathcal P_N$ by inequalities, consider a proper
cyclic interval $H\subsetneq E$. If $e\notin H$, the bound
$x(H)\le \overline b(H)$ becomes $y(H)\le \overline b(H)$. If $e\in H$, its
complement $E\setminus H$ is a nonempty interval in
$N$, and the same bound becomes
$y(E\setminus H)\ge r-\overline b(H)=\underline b(E\setminus H)$.
As $H$ ranges over the proper cyclic intervals of $E$, these are
exactly the upper and lower bounds on the intervals
of $N$. For each nonempty interval $I\subseteq N$,
define
\begin{equation}\label{eq:1-2}
 \mathcal P_I\coloneqq\{y\in\R^I:
 \underline b(J)\le y(J)\le \overline b(J)\text{ for every nonempty interval }J\subseteq I\}.
\end{equation}
For $I=N$, this agrees with the projection in \eqref{eq:1-3}.

\begin{lemma}\label{lem:projection}
For every nonempty interval $I\subseteq N$,
$\mathcal P_I=\operatorname{pr}_I(\mathcal P)$. For every positive integer $t$, each integral point
of $t\mathcal P_I$ extends to an integral point of $t\mathcal P$.
\end{lemma}

\begin{proof}
Rotate and label the cyclic order so that $E=[n]$ and $I=[k]$,
where $k\coloneqq|I|<n$. For $x\in \mathcal P$, write
$z_j\coloneqq\sum_{i=1}^j x_i$ for $1\le j\le n$, and let $z_0=0$;
then $z_n=r$.
For each pair of integers $0\le u,v\le n$, define
\[
 d_{uv}\coloneqq\max_{x\in\mathcal P}(z_v-z_u).
\]
For $u<v$, write $(u,v]=\{u+1,\ldots,v\}$; then
$d_{uv}=\overline b((u,v])$ and $d_{vu}=-\underline b((u,v])$.
Also $d_{uu}=0$ for $0\le u\le n$, so every $d_{uv}$ is an integer.
The conditions $x(E)=r$ and $x(H)\le \overline b(H)$ for every cyclic
interval $H\subseteq E$ become
\[
 z_0=0,\qquad z_n=r,\qquad z_v-z_u\le d_{uv},
\]
where the inequalities hold for all integers $0\le u,v\le n$.
For all integers $0\le u,v,w\le n$ and every $x\in\mathcal P$,
$z_w-z_u=(z_v-z_u)+(z_w-z_v)\le d_{uv}+d_{vw}$.
Taking the maximum over $x\in\mathcal P$ gives
$d_{uw}\le d_{uv}+d_{vw}$.

Fix a positive integer $t$ and $y\in t\mathcal P_I$.
Prescribe $z_0=0$, $z_a=\sum_{i=1}^a y_i$ for $1\le a\le k$,
and $z_n=rt$. Let $F\coloneqq\{0,\ldots,k,n\}$.
We claim that $z_v-z_u\le td_{uv}$ for all $u,v\in F$.
When $u=v$, equality holds because $d_{uu}=0$.
For $0\le u<v\le k$, the bounds defining $t\mathcal P_I$ give
\[
\begin{aligned}
 z_v-z_u&=y((u,v])\le t\overline b((u,v])=td_{uv},\\
 z_u-z_v&=-y((u,v])\le-t\underline b((u,v])=td_{vu}.
\end{aligned}
\]
For $1\le a\le k$, the bound
$t\underline b([1,a])\le z_a\le t\overline b([1,a])$ follows from
$y\in t\mathcal P_I$. Since
$d_{an}=r-\underline b([1,a])$ and $d_{na}=\overline b([1,a])-r$, we obtain
\[
 z_n-z_a=rt-z_a\le td_{an},\qquad
 z_a-z_n=z_a-rt\le td_{na}.
\]
The remaining pair $0,n$ satisfies
$z_n-z_0=rt=td_{0n}$ and $z_0-z_n=-rt=td_{n0}$.
This proves the claim for all ordered pairs in $F$.

For each $v\in\{0,\ldots,n\}$, define
\[
 \widetilde z_v\coloneqq\min_{a\in F}(z_a+td_{av}).
\]
For $v\in F$, the verified inequalities give
$z_a+td_{av}\ge z_v$ for every $a\in F$, with equality at
$a=v$. Hence $\widetilde z_v=z_v$ on $F$.
For any $u,v\in\{0,\ldots,n\}$, the inequalities
$d_{av}\le d_{au}+d_{uv}$ imply
\[
 \widetilde z_v
 \le\min_{a\in F}(z_a+td_{au}+td_{uv})
 =\widetilde z_u+td_{uv}.
\]
Define $x\in\R^E$ by $x_i=\widetilde z_i-\widetilde z_{i-1}$
for $1\le i\le n$. Since $\widetilde z=z$ on $F$, we have
$x_i=y_i$ for $i\le k$ and $x(E)=rt$.
The partial sums $\widetilde z_i/t$ of $x/t$ satisfy all the
difference inequalities defining $\mathcal P$, so $x\in t\mathcal P$.
At $t=1$, this gives $\mathcal P_I\subseteq\operatorname{pr}_I(\mathcal P)$.
Conversely, the tight bounds on $\mathcal P$ give
$\underline b(J)\le x(J)\le \overline b(J)$ for every $x\in\mathcal P$ and
every interval $J\subseteq I$, so $\operatorname{pr}_I(\mathcal P)\subseteq\mathcal P_I$.
Thus $\mathcal P_I=\operatorname{pr}_I(\mathcal P)$.

If $y$ is integral, then every prescribed $z_a$ and every
$td_{av}$ is an integer. The defining minimum gives
$\widetilde z_v\in\Z$ for all $v$, and consequently
$x\in t\mathcal P\cap\Z^E$.
\end{proof}

Coordinate projection preserves the sum on $I$. Thus
Lemma~\ref{lem:projection} and the definitions in \eqref{eq:1-1}
show that the minimum and maximum of $y(I)$ on $\mathcal P_I$
are $\underline b(I)$ and $\overline b(I)$, respectively.
For $I=N$, the identity $x(N)=r-x_e$ and the range $0\le x_e\le c_e$
give
\[
 \min_{y\in \mathcal P_N}y(N)=r-c_e,\qquad
 \max_{y\in \mathcal P_N}y(N)=r.
\]

\begin{definition}\label{def:interval-generating-function}
For a nonempty interval $I\subseteq N$ and a positive
integer $t$, define the \emph{lattice-point generating function}
in the indeterminate $z$ by
\begin{equation}\label{eq:interval-generating-function}
 \beta_I(t,z)\coloneqq\sum_{y\in t\mathcal P_I\cap\Z^I}z^{y(I)}.
\end{equation}
The exponent records the coordinate sum $y(I)$.
Each projected lattice point $y$ contributes once, even when
several points of $t\mathcal P\cap\Z^E$ project to $y$.
\end{definition}

For a polynomial or formal power series $\varphi(z)$, we write
$[z^k]\varphi(z)$ for the coefficient of $z^k$.

For $i\in N$ and a positive integer $t$, the integer points of
$t\mathcal P_{\{i\}}=[0,c_it]$ are $0,1,\ldots,c_it$. Hence
\begin{equation}\label{eq:singleton-generating-function}
 \beta_{\{i\}}(t,z)=\sum_{a=0}^{c_it}z^a
 =1+z+\cdots+z^{c_it},
\end{equation}
where the exponent records the value of coordinate $i$. For the
deleted coordinate $e$, write $\beta_{\{e\}}(t,z)\coloneqq\sum_{a=0}^{c_e t}z^a$.

\begin{lemma}\label{lem:root-count}
For every positive integer $t$,
\begin{equation}\label{eq:root-count}
 i(\mathcal P,t)
 =[z^{rt}]\bigl(\beta_N(t,z)\beta_{\{e\}}(t,z)\bigr).
\end{equation}
\end{lemma}

\begin{proof}
Expanding the product on the right-hand side of \eqref{eq:root-count}
shows that the coefficient of $z^{rt}$ counts pairs $(y,a)$ with
$y\in t\mathcal P_N\cap\Z^N$, $a\in\{0,\ldots,c_et\}$, and
$y(N)+a=rt$.
Since $\underline b(N)=r-c_e$ and $\overline b(N)=r$, every
$y\in t\mathcal P_N\cap\Z^N$ satisfies
\[
 (r-c_e)t\le y(N)\le rt.
\]
Thus $a=rt-y(N)$ is the unique integer in $\{0,\ldots,c_et\}$
for which $y(N)+a=rt$.
By \eqref{eq:1-3} after dilation by $t$, adjoining the coordinate
$x_e=rt-y(N)$ gives a bijection from $t\mathcal P_N\cap\Z^N$
to $t\mathcal P\cap\Z^E$, with inverse $\operatorname{pr}_N$.
Consequently,
\[
 [z^{rt}]\bigl(\beta_N(t,z)\beta_{\{e\}}(t,z)\bigr)
 =\#(t\mathcal P_N\cap\Z^N)
 =\#(t\mathcal P\cap\Z^E)
 =i(\mathcal P,t).
\]
\end{proof}

\begin{example}
\label{ex:running-polypositroid}
All the examples below use the polytope
\[
 \mathcal P_\star
 =\{x\in[0,2]^2\times[0,1]^3:
       x_1+\cdots+x_5=3,\ x_1+x_2+x_3\le2\},
\]
with ground set $E=\{1,\ldots,5\}$ in the cyclic order $1,2,3,4,5$.
To verify that $\mathcal P_\star$ is a generalized permutohedron,
define, for each $A\subseteq E$,
\[
 \overline b_\star(A)=\min\bigl\{3,\,
   \min\{2,\,2|A\cap\{1,2\}|+|A\cap\{3\}|\}
      +|A\cap\{4,5\}|\bigr\}.
\]
The functions $A\mapsto2|A\cap\{1,2\}|+|A\cap\{3\}|$ and
$A\mapsto|A\cap\{4,5\}|$ are nondecreasing, integer-valued, and modular.
Sums and truncations at nonnegative integer constants preserve
nondecreasing submodularity and integer values, so $\overline b_\star$ is
an integer-valued submodular function with
$\overline b_\star(\varnothing)=0$ and $\overline b_\star(E)=3$.

We claim that $\mathcal P_\star$ is the base polytope of $\overline b_\star$.
For $x\in\mathcal P_\star$, nonnegativity, the coordinate upper bounds,
$x(\{1,2,3\})\le2$, and $x(E)=3$ give
$x(A)\le \overline b_\star(A)$ for every $A\subseteq E$.
Conversely, let $x\in\R^5$ satisfy $x(E)=3$ and
$x(A)\le \overline b_\star(A)$ for every $A\subseteq E$.
Taking $A=\{i\}$ gives $x_1,x_2\le2$ and $x_3,x_4,x_5\le1$.
For each $i\in E$, taking $A=E\setminus\{i\}$ gives
\[
 x_i=3-x(E\setminus\{i\})
 \ge3-\overline b_\star(E\setminus\{i\})=0.
\]
Taking $A=\{1,2,3\}$ gives $x_1+x_2+x_3\le2$, so $x\in\mathcal P_\star$.
By the characterization following
\eqref{eq:full-support}, $\mathcal P_\star$ is an integral generalized
permutohedron.
Under $x(E)=3$, each lower bound $x_i\ge0$ is equivalent to
$x(E\setminus\{i\})\le3$. The sets $\{i\}$,
$E\setminus\{i\}$, and $\{1,2,3\}$ are cyclic intervals,
so $\mathcal P_\star$ is an integral polypositroid.

In the examples, take $\mathcal P=\mathcal P_\star$, $e=5$, and
$N=\{1,2,3,4\}$. We retain the notations $\underline b,\overline b$
for the bounds of $\mathcal P_\star$ and $\mathcal P_I,\beta_I$ for its
projection and generating function on each interval $I\subseteq N$.
The coordinate maxima are
$c_1=c_2=2$ and $c_3=c_4=c_5=1$.
Deleting coordinate $5$ gives
\[
 \mathcal P_N=\{y\in[0,2]^2\times[0,1]^2:
       y_1+y_2+y_3\le2,\ 2\le y(N)\le3\}.
\]
The inverse of $\operatorname{pr}_N$ adjoins $x_5=3-y(N)$, which lies in
$[0,1]$ precisely when $2\le y(N)\le3$.
For $I=\{1,2\}$, the projections onto $I$ and $\{3,4\}$ are
\[
 \mathcal P_I=\{(y_1,y_2)\in\R_{\ge0}^2:y_1+y_2\le2\},
 \qquad \mathcal P_{\{3,4\}}=[0,1]^2.
\]
For every positive integer $t$, Lemma~\ref{lem:root-count}
applied to $\mathcal P_\star$ with $r=3$, $e=5$, and $c_5=1$ gives
\[
 i(\mathcal P_\star,t)=[z^{3t}]\,\beta_N(t,z)(1+z+\cdots+z^t).
\]
For $x\in t\mathcal P_\star\cap\Z^5$, put $s=x_1+x_2+x_3$.
The inequalities $s\le2t$ and $3t-s=x_4+x_5\le2t$ give
$t\le s\le2t$.
For each integer $s$ in this range and each $a\in\{0,\ldots,t\}$,
setting $x_3=a$ leaves $s-a+1$ choices for $(x_1,x_2)$.
The bounds $x_1,x_2\le2t$ hold automatically because
$0\le x_1+x_2=s-a\le2t$.
The pairs $(x_4,x_5)$ with sum $3t-s$ and both entries in
$\{0,\ldots,t\}$ are obtained by choosing
$x_4\in\{2t-s,\ldots,t\}$, giving $s-t+1$ choices.
Consequently,
\[
\begin{aligned}
 i(\mathcal P_\star,t)
 &=\sum_{s=t}^{2t}(s-t+1)\sum_{a=0}^t(s-a+1)\\
 &=\frac{(t+1)^2(t+2)(7t+6)}{12}\\
 &=\frac7{12}t^4+\frac{17}{6}t^3+\frac{59}{12}t^2
   +\frac{11}{3}t+1.
\end{aligned}
\]
\end{example}

\section{Interval recursion and finite normal form}\label{sec:counting}

Throughout this section, $I$ denotes a nonempty interval
in the linear order on $N=E\setminus\{e\}$ fixed in Section~\ref{sec:geometry}.
We derive a recursion for $\beta_I(t,z)$ by cancelling terms
indexed by interval partitions. We then encode the generating
functions by polynomial sequences with finitely many nonzero
components and express the coefficient of $z^{ht}$ using these polynomials.

\subsection{An inclusion-exclusion recursion}\label{sec:violations}

An \emph{interval partition} of $I$ is a partition into nonempty
intervals. Write $\operatorname{Int}(I)$ for the set of
interval partitions of $I$.
For $\pi\in\operatorname{Int}(I)$, the product $\prod_{J\in\pi}\mathcal P_J$
imposes the bounds on intervals contained in a block of $\pi$.
A vector in the product can still violate a bound on an interval
meeting two blocks.

Fix $y=(y_i)_{i\in I}\in\R^I$, without assuming
$y\in\mathcal P_I$.
\begin{definition}\label{def:bad-intervals}
Relative to $y$, a nonempty interval $J\subseteq I$ is
\emph{upper-violating} if $y(J)>\overline b(J)$, and
\emph{lower-violating} if $y(J)<\underline b(J)$. We call $J$
\emph{violating} if $y(J)\notin[\underline b(J),\overline b(J)]$.
A violating interval $J$ is \emph{minimal} if every nonempty
proper subinterval of $J$ is not violating.

Order the minimal violating intervals $J_1,\ldots,J_m$ by increasing
left endpoint. A \emph{group} is a maximal nonempty consecutive sublist
$J_a,\ldots,J_d$ such that $J_k\cap J_{k+1}\ne\varnothing$ for every
$a\le k<d$. The union of each group is an interval. Since no minimal
violating interval contains another, both endpoints increase strictly
along the list, so the unions of distinct groups are disjoint.
\end{definition}

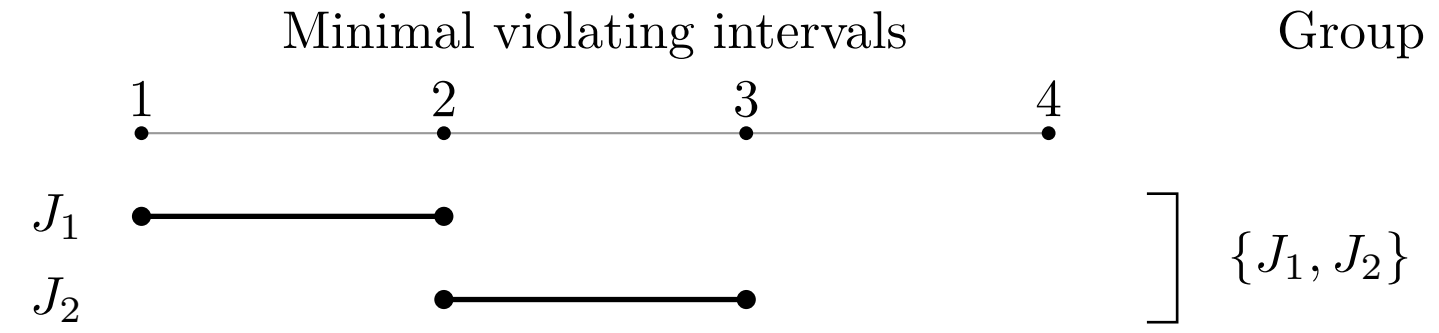
\begin{figure}[htbp]
\centering
\begin{tikzpicture}[x=1cm,y=1cm,font=\small]
 \node at (3,1.5) {Minimal violating intervals};
 \draw[black!40,line width=0.4pt] (0,0.85)--(6,0.85);
 \foreach \i/\x in {1/0,2/2,3/4,4/6} {
   \fill (\x,0.85) circle (1.3pt);
   \node[above,inner sep=3pt] at (\x,0.85) {$\i$};
 }
 \node[anchor=east] at (-0.25,0.3) {$J_1$};
 \draw[line width=1pt] (0,0.3)--(2,0.3);
 \fill (0,0.3) circle (1.8pt);
 \fill (2,0.3) circle (1.8pt);
 \node[anchor=east] at (-0.25,-0.25) {$J_2$};
 \draw[line width=1pt] (2,-0.25)--(4,-0.25);
 \fill (2,-0.25) circle (1.8pt);
 \fill (4,-0.25) circle (1.8pt);
 \node at (8,1.5) {Group};
 \draw[line width=0.5pt]
   (6.65,0.45)--(6.85,0.45)--(6.85,-0.4)--(6.65,-0.4);
 \node[anchor=west] at (7.05,0.025) {$\{J_1,J_2\}$};
\end{tikzpicture}
\caption{The intersecting minimal violating intervals for
$y=(a,2,b,0)$ in Example~\ref{ex:running-overlap} form a single group.}
\label{fig:interval-groups}
\end{figure}

\begin{example}\label{ex:running-overlap}
For the polypositroid in Example~\ref{ex:running-polypositroid},
take $I=N$. For real numbers $a,b>0$ with $a+b\le1$, let
\[
 y=(a,2,b,0).
\]
All singletons satisfy their bounds, and
$y(\{3,4\})=b\in[0,2]$.
The intervals $J_1=\{1,2\}$ and $J_2=\{2,3\}$ are upper-violating,
since
\[
 y(J_1)=2+a>\overline b(J_1)=2,\qquad
 y(J_2)=2+b>\overline b(J_2)=2.
\]
Every interval of length at least three contains $J_1$ or $J_2$,
so $J_1,J_2$ are exactly the minimal violating intervals.
They intersect at $2$ and form one group, as shown in
Figure~\ref{fig:interval-groups}.
The total sum satisfies
$y(I)=2+a+b\in[\underline b(I),\overline b(I)]=[2,3]$.
\end{example}

Fix an interval partition $I=K_1\sqcup\cdots\sqcup K_q$, with $q\ge1$
and the blocks listed in the linear order on $I$.
The proper intervals formed from whole blocks constitute
\[
 \mathcal H\coloneqq
 \{K_a\cup\cdots\cup K_d:1\le a\le d\le q,\ (a,d)\ne(1,q)\}.
\]
We test $\underline b(H)\le y(H)\le\overline b(H)$ only for
$H\in\mathcal H$, and treat the bounds on $I$ separately.
A violating $H\in\mathcal H$ is \emph{minimal in $\mathcal H$}
if no $G\in\mathcal H$ with $G\subsetneq H$ is violating.
List the minimal violating members of $\mathcal H$ as $J_1,\ldots,J_m$
in increasing order of their left endpoints.

\begin{samepage}
\begin{lemma}\label{lem:minimal-violations}
Let $\mathcal J$ be one of the groups of $J_1,\ldots,J_m$
in the sense of Definition~\ref{def:bad-intervals}, and put
$U\coloneqq\bigcup_{J\in\mathcal J}J$.
Then one of the following holds.
\begin{enumerate}[label=(\roman*),leftmargin=*,topsep=3pt,itemsep=2pt,parsep=0pt]
\item $y(J)>\overline b(J)$ for every $J\in\mathcal J$, and $y(U)>\overline b(U)$.
\item $y(J)<\underline b(J)$ for every $J\in\mathcal J$, and $y(U)<\underline b(U)$.
\end{enumerate}
If $\underline b(I)\le y(I)\le \overline b(I)$, then $U\subsetneq I$.
\end{lemma}
\end{samepage}

\begin{proof}
Suppose that two intersecting minimal violating members $A,B\in\mathcal H$
satisfy $y(A)>\overline b(A)$ and $y(B)<\underline b(B)$.
Neither interval contains the other, since the larger would
contain a violating member of $\mathcal H$ as a proper subset.
Thus $A\setminus B$ and $B\setminus A$ belong to $\mathcal H$
and are proper subsets of $A$ and $B$, respectively.
By minimality of $A$ and $B$, additivity of coordinate sums,
and \eqref{eq:1-4}, we obtain
\[
\begin{aligned}
 \overline b(A)-\underline b(B)&<y(A)-y(B)\\
 &=y(A\setminus B)-y(B\setminus A)\\
 &\le \overline b(A\setminus B)-\underline b(B\setminus A)\\
 &\le \overline b(A)-\underline b(B),
\end{aligned}
\]
a contradiction. Successive members of $\mathcal J$ intersect,
so all members are upper-violating or all are lower-violating.

Write $\mathcal J=\{J_a,\ldots,J_d\}$, with $1\le a\le d\le m$.
For $a<k\le d$, put
$H\coloneqq J_a\cup\cdots\cup J_{k-1}$ and $J\coloneqq J_k$.
Successive members of $\mathcal J$ intersect, so $H$ is a union
of consecutive $K_i$. Since the left and right endpoints of
$J_a,\ldots,J_d$ increase strictly, we have
\[
 \varnothing\ne H\cap J=J_{k-1}\cap J_k\subsetneq J,
\]
where strict containment holds because the right endpoint of
$J_{k-1}$ precedes the right endpoint of $J_k$.
The intersection $H\cap J$ belongs to $\mathcal H$, so
minimality of $J$ in $\mathcal H$ gives
\[
 \underline b(H\cap J)\le y(H\cap J)\le \overline b(H\cap J).
\]

If all members of $\mathcal J$ are upper-violating, start with
$y(J_a)>\overline b(J_a)$ and assume $y(H)>\overline b(H)$ at the induction step.
Using $y(J)>\overline b(J)$, the upper bound
$y(H\cap J)\le \overline b(H\cap J)$, and submodularity of $\overline b$, we obtain
\begin{equation}\label{eq:3-1}
\begin{aligned}
 y(H\cup J)&=y(H)+y(J)-y(H\cap J)\\
 &>\overline b(H)+\overline b(J)-\overline b(H\cap J)\ge \overline b(H\cup J).
\end{aligned}
\end{equation}
Induction on $k$ yields $y(U)>\overline b(U)$.
If all members of $\mathcal J$ are lower-violating, the induction
based on \eqref{eq:3-1}, applied to $-y$ and the submodular function
$-\underline b$, yields $-y(U)>-\underline b(U)$, hence
$y(U)<\underline b(U)$.
Since $U$ is violating, it cannot equal $I$ when
$\underline b(I)\le y(I)\le\overline b(I)$; thus $U\subsetneq I$.
\end{proof}

Our goal is to compute $\beta_I(t,z)$ from the generating functions of
shorter intervals by inclusion--exclusion. After imposing the bounds on
the total coordinate sum, the signed sum over partitions must give
multiplicity one to each lattice point of $t\mathcal P_I$ and zero to
every other lattice vector.

Assume $q\ge2$ and number the gaps between
$K_1,\ldots,K_q$ from left to right by $1,\ldots,q-1$.
For $C\subseteq\{1,\ldots,q-1\}$, let $\pi(C)$ be the partition of $I$
whose block boundaries are precisely the gaps indexed by $C$.
The map $C\mapsto\pi(C)$ is a bijection from subsets of the gaps to
partitions obtained by merging consecutive $K_i$, and $|\pi(C)|=|C|+1$.
An interval $H=K_a\cup\cdots\cup K_d\in\mathcal H$ is contained in one block
of $\pi(C)$ if and only if $C\cap\{a,\ldots,d-1\}=\varnothing$.

Let $\mathcal C_y$ be the family of nonempty subsets
$C\subseteq\{1,\ldots,q-1\}$ satisfying
\[
 C\cap\{a,\ldots,d-1\}\ne\varnothing
\]
for every $H=K_a\cup\cdots\cup K_d$ minimal among the violating members
of $\mathcal H$.
Every violating member of $\mathcal H$ contains a minimal one, so a gap
inside the smaller interval also lies inside the larger interval.
If a single block $K_i$ is violating, then $\mathcal C_y=\varnothing$,
because $\{a,\ldots,d-1\}$ is empty when $a=d$.
The \emph{signed multiplicity} of $y$, relative to the fixed
partition $K_1,\ldots,K_q$, is
\begin{equation}\label{eq:4-3}
 \sigma(y)\coloneqq\sum_{C\in\mathcal C_y}(-1)^{|C|+1}.
\end{equation}

\begin{lemma}\label{lem:cut-sum}
If $\underline b(I)\le y(I)\le\overline b(I)$, then
\[
 \sigma(y)=
 \begin{cases}
  1,&\text{if no member of }\mathcal H\text{ is violating},\\
  0,&\text{if at least one member of }\mathcal H\text{ is violating}.
 \end{cases}
\]
Moreover, for every $y\in\R^I$,
\begin{equation}\label{eq:geometric-recursion}
 \mathbf1_{\mathcal P_I}(y)
 =\mathbf1_{[\underline b(I),\overline b(I)]}(y(I))
  \sum_{\substack{\pi\in\operatorname{Int}(I)\\|\pi|\ge2}}
   (-1)^{|\pi|}\mathbf1_{\prod_{J\in\pi}\mathcal P_J}(y).
\end{equation}
\end{lemma}

\begin{proof}
If no member of $\mathcal H$ is violating, every nonempty subset
of $\{1,\ldots,q-1\}$ belongs to $\mathcal C_y$, so
\[
 \sigma(y)=\sum_{j=1}^{q-1}(-1)^{j+1}\binom{q-1}{j}=1.
\]
Suppose instead that $\mathcal H$ has a violating member. Since
$\underline b(I)\le y(I)\le\overline b(I)$, Lemma~\ref{lem:minimal-violations}
implies that the union of each group of minimal violating members is a
proper subset of $I$.
Choose one group $\mathcal J$ and put
$U\coloneqq\bigcup_{J\in\mathcal J}J$.
Choose a gap $g$ between a block of $U$ and an adjacent block
of $I\setminus U$. The unions of different groups are disjoint,
so $g$ lies inside no minimal violating member of $\mathcal H$.
Define
\[
 \iota_g:\mathcal C_y\longrightarrow\mathcal C_y,\qquad
 \iota_g(C)\coloneqq C\mathbin{\triangle}\{g\}.
\]
Adding or removing $g$ preserves the condition defining
$\mathcal C_y$, since it does not change whether $C$ cuts
any minimal violating member of $\mathcal H$.
For each $C\in\mathcal C_y$, at least one gap other than $g$
is needed to cut a minimal violating member. Thus $\iota_g(C)$
is nonempty. We have $\iota_g(\iota_g(C))=C$, and
$|\iota_g(C)|$ differs from $|C|$ by one. The terms indexed
by $C$ and $\iota_g(C)$ in \eqref{eq:4-3} therefore have opposite
signs, giving $\sigma(y)=0$.
The binomial sum in the first case requires no bounds on $y(I)$.
In the second case, whenever every group union is a proper subset of $I$,
we can choose the gap $g$ and apply $\iota_g$ without those bounds.
Thus for every $y\in\R^I$,
\begin{equation}\label{eq:4-4}
 \sigma(y)<0\quad\Longrightarrow\quad
 \exists\,\mathcal J:\quad\bigcup_{J\in\mathcal J}J=I,
\end{equation}
where $\mathcal J$ is a group of minimal violating members of $\mathcal H$.

To prove \eqref{eq:geometric-recursion}, take the partition of
$I$ into singleton blocks. Then $\mathcal H$ consists of all
nonempty proper intervals of $I$.
For a nonempty subset $C$ of the gaps between consecutive
elements of $I$, the product $\prod_{J\in\pi(C)}\mathcal P_J$
contains $y$ exactly when every interval contained in a block
of $\pi(C)$ satisfies its bounds, or equivalently, when
$C\in\mathcal C_y$. The bijection $C\mapsto\pi(C)$ and
$|\pi(C)|=|C|+1$ identify the partition sum in
\eqref{eq:geometric-recursion} with $\sigma(y)$ for the singleton
partition. If $\underline b(I)\le y(I)\le\overline b(I)$,
the values of $\sigma(y)$ just proved give
\eqref{eq:geometric-recursion}. If $y(I)$ lies outside this
closed interval, both sides of \eqref{eq:geometric-recursion}
vanish.
\end{proof}

In \eqref{eq:4-4}, every member of $\mathcal J$ is a proper
subset of $I$, so $|\mathcal J|\ge2$.
Lemma~\ref{lem:minimal-violations} also shows that the members
of $\mathcal J$ are all upper-violating or all lower-violating.

\begin{example}\label{ex:running-cuts}
Use $y=(a,2,b,0)$ from Example~\ref{ex:running-overlap}, with
$a,b>0$ and $a+b\le1$, and take the four singletons as fixed blocks.
The minimal violating intervals $\{1,2\}$ and $\{2,3\}$ force
every $C\in\mathcal C_y$ to contain gaps $1$ and $2$.
Gap $3$, between $\{3\}$ and $\{4\}$, lies inside neither
minimal violating interval, so
\[
 \mathcal C_y=\{\{1,2\},\{1,2,3\}\}.
\]
The involution $\iota_3$ pairs the corresponding partitions:
\[
 \underbrace{\{1\}\mid\{2\}\mid\{3,4\}}_{
   \substack{C=\{1,2\}\\(-1)^{|\pi(C)|}=-1}}
 \quad\overset{\iota_3}{\longleftrightarrow}\quad
 \underbrace{\{1\}\mid\{2\}\mid\{3\}\mid\{4\}}_{
   \substack{C=\{1,2,3\}\\(-1)^{|\pi(C)|}=+1}}.
\]
Hence $\sigma(y)=-1+1=0$. Since $y(I)$ satisfies its bounds, the
singleton-block case of Lemma~\ref{lem:cut-sum} gives
$\mathbf1_{\mathcal P_I}(y)=0$, so $y\notin\mathcal P_I$.
The endpoint equalities $y_2=\overline b(\{2\})=2$ and
$y_4=\underline b(\{4\})=0$ do not cause violations.
The same pair of partitions remains admissible when $a+b=1$,
so that $y(I)=\overline b(I)=3$.
\end{example}

We now pass from \eqref{eq:geometric-recursion} to generating
functions. For a polynomial $\varphi(z)$ and integers $0\le a\le d$,
we retain only the terms whose exponents lie in $[a,d]$ by setting
\[
 \left.\varphi(z)\right|_{[a,d]}
 \coloneqq\sum_{k=a}^{d}([z^k]\varphi(z))z^k.
\]
For $|I|\ge2$, let
\[
 \gamma_I(t,z)\coloneqq
 \sum_{\substack{\pi\in\operatorname{Int}(I)\\|\pi|\ge2}}
 (-1)^{|\pi|}\prod_{J\in\pi}\beta_J(t,z).
\]
Thus $\gamma_I(t,z)$ is the signed sum of the block products
before imposing the bounds on the total coordinate sum over $I$.

\begin{proposition}\label{prop:exact-recursion}
If $|I|\ge2$, then for every positive integer $t$,
\begin{equation}\label{eq:counting-recursion}
 \beta_I(t,z)=
 \left.\gamma_I(t,z)\right|_{[t\underline b(I),\,t\overline b(I)]}.
\end{equation}
\end{proposition}

\begin{proof}
Fix a positive integer $t$ and a partition $\pi\in\operatorname{Int}(I)$
with $|\pi|\ge2$. Choosing a lattice point of $t\mathcal P_J$
independently for each $J\in\pi$ determines a unique $a\in\Z_{\ge0}^I$ whose
projection $\operatorname{pr}_J(a)$ is the chosen point on $J$.
Since $a(I)=\sum_{J\in\pi}a(J)$, we obtain
\[
\begin{aligned}
 \prod_{J\in\pi}\beta_J(t,z)
 &=\prod_{J\in\pi}\left(\sum_{x\in t\mathcal P_J\cap\Z^J}z^{x(J)}\right)\\
 &=\sum_{a\in\prod_{J\in\pi}(t\mathcal P_J\cap\Z^J)}z^{a(I)}.
\end{aligned}
\]
For every $a\in\Z_{\ge0}^I$, substituting $y=a/t$ in
\eqref{eq:geometric-recursion} gives
\[
 \mathbf1_{t\mathcal P_I}(a)
 =\mathbf1_{[t\underline b(I),\,t\overline b(I)]}(a(I))
   \sum_{\substack{\pi\in\operatorname{Int}(I)\\|\pi|\ge2}}
   (-1)^{|\pi|}
   \prod_{J\in\pi}\mathbf1_{t\mathcal P_J}(\operatorname{pr}_J(a)).
\]
Multiply by $z^{a(I)}$ and sum over $a\in\Z_{\ge0}^I$.
Since $\mathcal P_I\subseteq\R_{\ge0}^I$, the left-hand side gives
\[
 \sum_{a\in\Z_{\ge0}^I}\mathbf1_{t\mathcal P_I}(a)z^{a(I)}
 =\sum_{a\in t\mathcal P_I\cap\Z^I}z^{a(I)}
 =\beta_I(t,z).
\]
On the right, the signed sum of the block products is $\gamma_I(t,z)$,
and the interval indicator retains exactly its terms with exponents in
$[t\underline b(I),\,t\overline b(I)]$. This gives
\eqref{eq:counting-recursion}.
\end{proof}

The singleton polynomials $\beta_{\{i\}}(t,z)$, for $i\in N$,
are the initial values for \eqref{eq:counting-recursion}, since
every block $J\in\pi$ defining $\gamma_I$ satisfies $|J|<|I|$.

\subsection{Finite normal forms and coefficient extraction}\label{subsec:normal-forms}

We use polynomial sequences with finitely many nonzero components
to represent the products and signed sums defining $\gamma_I(t,z)$.
For positive integers $c$ and $t$, the generating polynomial for the integers
$a$ satisfying $0\le a\le ct$ is
\[
 \beta_{[0,c]}(t,z)\coloneqq\sum_{a=0}^{ct}z^a=\frac{1-z^{ct+1}}{1-z}.
\]
By \eqref{eq:singleton-generating-function},
$\beta_{\{i\}}(t,z)=\beta_{[0,c_i]}(t,z)$ for each $i\in N$.
A single nonnegative integer coordinate contributes one term $z^a$ for
each possible value $a$, so its generating series is
$x=\sum_{a\ge0}z^a=(1-z)^{-1}$. To impose $a\le ct$, we subtract
the tail $z^{ct+1}x$ from $x$. Using $zx=x-1$, we obtain
\begin{equation}\label{eq:beta-c-expansion}
 \beta_{[0,c]}(t,z)=x-z^{ct+1}x=x+z^{ct}(1-x).
\end{equation}

For each positive integer $t$, let $\beta(t,z)$ be a formal power
series in $z$ with real coefficients.
Let $x$ be an independent indeterminate.

\begin{definition}\label{def:finite-normal-form}
A \emph{finite normal form} of $\beta(t,z)$ is a sequence of
polynomials in $\R[t,x]$,
\[
 \boldsymbol{\varphi}=(\varphi_j(t,x))_{j\ge0},
\]
with $\varphi_j=0$ for all but finitely many $j$, such that, for every
positive integer $t$, the following identity holds as a formal
power series in $z$:
\begin{equation}\label{eq:2-1}
 \beta(t,z)=\sum_{j\ge0}z^{jt}\varphi_j\left(t,\frac1{1-z}\right).
\end{equation}
\end{definition}

A summand in \eqref{eq:2-1} can have infinitely many nonzero
coefficients in $z$, some of them negative.

For example, take $\beta=\beta_{[0,c]}$.
Comparing \eqref{eq:beta-c-expansion} with \eqref{eq:2-1} gives a finite
normal form $\boldsymbol{\alpha}_c$ whose only nonzero components are
$(\boldsymbol{\alpha}_c)_0=x$ and $(\boldsymbol{\alpha}_c)_c=1-x$.
Expanding $x$ and $z^{ct}(1-x)$ in powers of $z$ gives
\[
 \beta_{[0,c]}(t,z)
 =\sum_{a=0}^{+\infty}z^a-\sum_{a=ct+1}^{+\infty}z^a
 =\sum_{a=0}^{ct}z^a.
\]

Let $\boldsymbol{\varphi}=(\varphi_j)_{j\ge0}$ and
$\widetilde{\boldsymbol{\varphi}}=(\widetilde\varphi_j)_{j\ge0}$ be finite
normal forms. For each integer $h\ge0$, define the $h$th components
of their sum and convolution product by
\[
\begin{aligned}
 (\boldsymbol{\varphi}+\widetilde{\boldsymbol{\varphi}})_h
   &\coloneqq \varphi_h+\widetilde\varphi_h,\\
 (\boldsymbol{\varphi}\widetilde{\boldsymbol{\varphi}})_h
   &\coloneqq \sum_{j=0}^h\varphi_j\widetilde\varphi_{h-j}.
\end{aligned}
\]
We claim that if $\boldsymbol{\varphi}$ and
$\widetilde{\boldsymbol{\varphi}}$ are finite normal forms of
$\beta$ and $\widetilde\beta$, respectively, then their convolution product
$\boldsymbol{\varphi}\widetilde{\boldsymbol{\varphi}}$ is a finite
normal form of $\beta\widetilde\beta$.
For each positive integer $t$, put $x=(1-z)^{-1}$.
Substituting the two normal forms and grouping the finitely many
nonzero terms by $h=a+b$ gives
\[
\begin{aligned}
 \beta(t,z)\widetilde\beta(t,z)
 &=\sum_{a,b\ge0}z^{(a+b)t}\varphi_a(t,x)\widetilde\varphi_b(t,x)\\
 &=\sum_{h\ge0}z^{ht}
   (\boldsymbol{\varphi}\widetilde{\boldsymbol{\varphi}})_h(t,x).
\end{aligned}
\]
Adding or subtracting the expansions of $\beta$ and $\widetilde\beta$
in \eqref{eq:2-1} shows that
$\boldsymbol{\varphi}\pm\widetilde{\boldsymbol{\varphi}}$ is a
finite normal form of $\beta\pm\widetilde\beta$.

In Section~\ref{subsec:positive-states}, we will use convolution,
componentwise addition and subtraction, and closed truncation to
construct finite normal forms of the interval generating functions
$\beta_I(t,z)$.

For each integer $k\ge0$, let $\binom{v}{k}$ denote the unique polynomial
in $v$ whose value at every integer $n\ge k$ is
$\frac{n!}{k!(n-k)!}$. In particular, $\binom{v}{0}=1$.

For an integer $b\ge1$, write $(1-z)^{-b}$ as a product of $b$ copies
of the single-coordinate series $\sum_{v\ge0}z^v$. Choosing $z^{v_i}$
from factor $i$ gives the monomial $z^{v_1+\cdots+v_b}$. For each
integer $s\ge0$, the coefficient of $z^s$ therefore counts the tuples
$(v_1,\ldots,v_b)\in\Z_{\ge0}^b$ satisfying $v_1+\cdots+v_b=s$.
By stars and bars,
\begin{equation}\label{eq:negative-binomial-coefficient}
 [z^s](1-z)^{-b}=\binom{s+b-1}{b-1}.
\end{equation}

\begin{lemma}\label{lem:coefficient-polynomial}
Suppose that $\boldsymbol{\varphi}$ is a finite normal form of
$\beta(t,z)$, and write
$\varphi_j(t,x)=\sum_{b\ge0}\alpha_{jb}(t)x^b$ with
$\alpha_{jb}\in\R[t]$.
For each integer $h\ge0$, the polynomial
\begin{equation}\label{eq:2-3}
 \zeta_h(t)\coloneqq\varphi_h(t,1)
   +\sum_{\substack{0\le j<h\\b\ge1}}\alpha_{jb}(t)
       \binom{(h-j)t+b-1}{b-1}
\end{equation}
is the unique polynomial in $\R[t]$ satisfying
$\zeta_h(t)=[z^{ht}]\beta(t,z)$
for every positive integer $t$.
\end{lemma}

\begin{proof}
The right-hand side of \eqref{eq:2-3} is a polynomial in $t$,
since only finitely many $\alpha_{jb}$ are nonzero and each
binomial factor is a polynomial in $t$.
Fix a positive integer $t$ and extract the coefficient of $z^{ht}$
in \eqref{eq:2-1}. Terms with $j>h$ contain only powers of degree at least $jt>ht$.
For $j=h$, the constant term of $(1-z)^{-1}$ is $1$, so the coefficient
is $\sum_{b\ge0}\alpha_{hb}(t)=\varphi_h(t,1)$. For $j<h$,
the required exponent $(h-j)t$ in $(1-z)^{-b}$ is positive, so
only terms with $b\ge1$ can occur. Applying \eqref{eq:negative-binomial-coefficient}
with $s=(h-j)t$ gives the binomial factors in \eqref{eq:2-3}, and hence
$[z^{ht}]\beta(t,z)=\zeta_h(t)$.
Any polynomial taking the values $[z^{ht}]\beta(t,z)$ at all positive
integers agrees with $\zeta_h$ at infinitely many points and hence equals $\zeta_h$.
\end{proof}

Every point of $t\mathcal P$ has coordinate sum $rt$, so we will take
$h=r$ when applying Lemma~\ref{lem:coefficient-polynomial} to a finite
normal form of $\beta(t,z)=\beta_N(t,z)\beta_{\{e\}}(t,z)$.
By \eqref{eq:root-count}, the Ehrhart polynomial $i(\mathcal P,t)$
has the values $[z^{rt}]\beta(t,z)$ at all positive integers $t$.
Uniqueness in Lemma~\ref{lem:coefficient-polynomial} therefore gives
\[
 i(\mathcal P,t)=\zeta_r(t).
\]
To prove that $\mathcal P$ is Ehrhart positive, it suffices to show
that $\zeta_r(t)$ has nonnegative coefficients.

\section{Integral operators and boundary limits}\label{sec:operators}

We represent finite normal forms by operators on functions of a
real variable. Divided differences give the integral kernels,
and limits of those kernels recover the corrections required
by closed truncation. We then construct the operator representing
$\beta_I(t,z)$ from the interval recurrence \eqref{eq:counting-recursion},
using operator composition for the products over blocks and closed
truncation for the bounds on the coordinate sum over $I$.

\subsection{Operators associated with finite normal forms}\label{subsec:integral-operators}

All integrals in this section are Lebesgue integrals.

\begin{definition}\label{def:auxiliary-function}
An \emph{auxiliary function} is a $1$-periodic function $\lambda:\R\to[0,1]$
satisfying the following two conditions.
\begin{enumerate}[label=(\roman*),leftmargin=*,topsep=4pt,itemsep=3pt,parsep=0pt]
\item The function $\lambda$ is \emph{piecewise affine} on one period.
There is a finite partition $0=a_0<a_1<\cdots<a_m=1$ such that
$\lambda$ is affine on each open interval $(a_{i-1},a_i)$.
\item The restriction of $\lambda$ to $[0,1]$ \emph{preserves Lebesgue measure}.
Equivalently, for every continuous function $\psi:[0,1]\to\R$,
\[
 \int_0^1\psi(\lambda(u))\,\mathrm{d}u=\int_0^1\psi(x)\,\mathrm{d}x.
\]
\end{enumerate}
\end{definition}

The fractional-part function $\lambda(u)=u-\lfloor u\rfloor$ is an
auxiliary function. Starting from this example or any other auxiliary
function $\lambda$, we obtain another by prescribing finitely many
point values per period. For distinct $u_1,\ldots,u_q\in[0,1)$ and
$x_1,\ldots,x_q\in[0,1]$, set $\widetilde\lambda(u_i+k)=x_i$
for $1\le i\le q$ and $k\in\Z$, and put
$\widetilde\lambda(u)=\lambda(u)$ elsewhere.

Let $\mathscr F$ be the vector space over $\R$ of locally bounded,
real-valued Borel measurable functions on $\R$ that vanish for all
sufficiently large arguments. An element $g$ of $\mathscr F[t]$ is a
polynomial in $t$ with coefficients in $\mathscr F$. We view $g$ as a map
$\R\to\R[t]$ by evaluating its coefficients at $u$ and leaving $t$
indeterminate. Thus
\[
 g(u)=\sum_{k=0}^d\vartheta_k(u)t^k,
 \qquad \vartheta_0,\ldots,\vartheta_d\in\mathscr F.
\]

Let $(D,\mu)$ be a measure space independent of $t$.
Let $\varphi(t,y)=\sum_{k=0}^d\alpha_k(y)t^k$ be a polynomial in $t$
with $\mu$-integrable coefficient functions $\alpha_k:D\to\R$.
Since the sum is finite, linearity of integration gives
\begin{equation}\label{eq:coefficientwise-integration}
 \int_D\varphi(t,y)\,\mathrm{d}\mu(y)
 =\sum_{k=0}^d\left(\int_D\alpha_k(y)\,\mathrm{d}\mu(y)\right)t^k.
\end{equation}

Given an auxiliary function $\lambda$, let $\mathsf S$ denote translation
by one and let $\mathsf X_\lambda$ act by multiplication by $\lambda$
plus $t$ times the tail integral. For $g\in\mathscr F[t]$ and $u\in\R$,
the two $\R[t]$-linear operators are given by
\[
\begin{aligned}
 \mathsf S,\mathsf X_\lambda &\vcentcolon \mathscr F[t]\longrightarrow\mathscr F[t],\\
 (\mathsf S g)(u) &\coloneqq g(u+1),\\
 (\mathsf X_\lambda g)(u) &\coloneqq \lambda(u)g(u)+t\int_u^{+\infty}g(v)\,\mathrm{d}v.
\end{aligned}
\]
Periodicity of $\lambda$ and the substitution $w=v+1$ in the tail
integral give $\mathsf S\mathsf X_\lambda=\mathsf X_\lambda\mathsf S$.

\begin{definition}\label{def:normal-form-operator}
Let $\lambda$ be an auxiliary function and let
$\boldsymbol{\varphi}=(\varphi_j)_{j\ge0}$ be a finite normal form.
Write $\varphi_j(t,x)=\sum_{b\ge0}\alpha_{jb}(t)x^b$ with
$\alpha_{jb}(t)\in\R[t]$.
In $\varphi_j(t,\mathsf X_\lambda)$ we substitute $\mathsf X_\lambda$ for $x$
and retain $t$ as an indeterminate. The $\R[t]$-linear operators
$\varphi_j(t,\mathsf X_\lambda)$ and $\mathsf R_{\lambda,\boldsymbol\varphi}$
from $\mathscr F[t]$ to $\mathscr F[t]$ are defined, for $g\in\mathscr F[t]$, by
\begin{equation}\label{eq:normal-form-operator}
\begin{aligned}
 \varphi_j(t,\mathsf X_\lambda)g
   &\coloneqq\sum_{b\ge0}\alpha_{jb}(t)\mathsf X_\lambda^{\,b}g,\\
 \mathsf R_{\lambda,\boldsymbol\varphi}g
   &\coloneqq\sum_{j\ge0}\mathsf S^j\bigl(\varphi_j(t,\mathsf X_\lambda)g\bigr).
\end{aligned}
\end{equation}
For $b\ge1$, $\mathsf X_\lambda^{\,b}$ denotes the $b$-fold composition
of $\mathsf X_\lambda$, and $\mathsf X_\lambda^0$ denotes the identity operator.
\end{definition}

\begin{lemma}\label{lem:normal-form-uniqueness}
Suppose that $\beta(t,z)$ has a finite normal form
$\boldsymbol{\varphi}=(\varphi_j)_{j\ge0}$ satisfying
\eqref{eq:2-1} for every positive integer $t$.
Then $\boldsymbol{\varphi}$ is uniquely determined by $\beta$.
\end{lemma}

\begin{proof}
Let $\boldsymbol\varphi=(\varphi_j)_{j\ge0}$ and
$\widetilde{\boldsymbol\varphi}=(\widetilde\varphi_j)_{j\ge0}$ be two finite
normal forms of $\beta(t,z)$.
Choose $d\in\Z_{\ge0}$ so that every nonzero polynomial among the
$\varphi_j$ and $\widetilde\varphi_j$ has degree at most $d$ in $x$.
For each $j\ge0$, set
\[
 \tau_j(t,z)\coloneqq(1-z)^d\left[
 \varphi_j\left(t,\frac1{1-z}\right)
 -\widetilde \varphi_j\left(t,\frac1{1-z}\right)\right].
\]
Each $\tau_j$ belongs to $\R[t,z]$ and has degree at most $d$ in $z$.
For each positive integer $n$, substitute $t=n$ in the two representations
\eqref{eq:2-1} of $\beta(t,z)$. Subtracting and multiplying by $(1-z)^d$ gives
the polynomial identity
\[
 \sum_{j\ge0}z^{jn}\tau_j(n,z)=0 \qquad\text{in }\R[z].
\]
For $n>d$, every exponent in $z^{jn}\tau_j(n,z)$ lies in $[jn,jn+d]$.
Since $jn+d<(j+1)n$, no other summand contains a power of $z$ in
$[jn,jn+d]$. Thus $\tau_j(n,z)$ is the zero polynomial for every $j$.
For each fixed $j$ and $0\le a\le d$, the polynomial $[z^a]\tau_j(t,z)$
has every integer $n>d$ as a root and is therefore the zero polynomial.
Hence each $\tau_j$ is the zero polynomial. Substituting $z=1-x^{-1}$
in the definition of $\tau_j$ and multiplying by $x^d$ gives
$\varphi_j(t,x)=\widetilde\varphi_j(t,x)$ for every $j$.
\end{proof}

Since $\beta$ determines $\boldsymbol\varphi$ by Lemma~\ref{lem:normal-form-uniqueness},
formula~\eqref{eq:normal-form-operator} determines
$\mathsf R_{\lambda,\boldsymbol\varphi}$ from $\beta$ and $\lambda$.

For an integer $c\ge1$, the finite normal form $\boldsymbol\alpha_c$ of
$\beta_{[0,c]}(t,z)$ from Section~\ref{subsec:normal-forms} has components
$x$ at index $0$ and $1-x$ at index $c$. For $c=c_i$, the corresponding operator gives the singleton case
$I=\{i\}$ of Proposition~\ref{prop:state-closure}. By \eqref{eq:normal-form-operator},
\[
 \mathsf R_{\lambda,\boldsymbol\alpha_c}
 =\mathsf X_\lambda+\mathsf S^c(1-\mathsf X_\lambda),
\]
where $1$ denotes the identity operator.
For $g\in\mathscr F[t]$ and $u\in\R$, the integrals over
$(u+c,+\infty)$ in
$(\mathsf X_\lambda g)(u)-(\mathsf X_\lambda g)(u+c)$
cancel. Using $\lambda(u+c)=\lambda(u)$ and setting $y=v-u$ gives
\begin{equation}\label{eq:2-13}
 (\mathsf R_{\lambda,\boldsymbol\alpha_c}g)(u)
 =\lambda(u)g(u)+(1-\lambda(u))g(u+c)
       +t\int_0^c g(u+y)\,\mathrm{d}y.
\end{equation}

\subsection{Divided differences and integral kernels}\label{subsec:divided-differences}

For a finite normal form $\boldsymbol{\varphi}=(\varphi_j)_{j\ge0}$
and an auxiliary function $\lambda$, we write
\[
 (\mathsf R_{\lambda,\boldsymbol\varphi} g)(u)
 =\sum_{j\ge0}\varphi_j(t,\lambda(u))g(u+j)
  +\int_u^{+\infty} K_{\boldsymbol{\varphi},\lambda}(u,v)g(v)\,\mathrm{d}v
\]
for every $g\in\mathscr F[t]$ and $u\in\R$. Using divided differences,
we derive the explicit kernel \eqref{eq:2-14} for a single operator
and the kernel formula \eqref{eq:2-16} for the composition of two such operators.

\begin{definition}\label{def:divided-differences}
For $\varphi(t,x)\in\R[t,x]$ and real \emph{nodes} $a_0,\ldots,a_k$,
let $\varphi[a_0,\ldots,a_k]$ denote the $k$th \emph{divided difference}
of $\varphi$ with respect to $x$ at these nodes. For $k=0$, set
$\varphi[a_0]\coloneqq\varphi(t,a_0)$.
For $k\ge1$ and $a_0\ne a_k$, the divided differences are
defined recursively by
\[
 \varphi[a_0,\ldots,a_k]
 \coloneqq\frac{\varphi[a_1,\ldots,a_k]-\varphi[a_0,\ldots,a_{k-1}]}{a_k-a_0}.
\]
For example, at the two nodes $0$ and $1$ we have
$\varphi[0,1]=\varphi(t,1)-\varphi(t,0)$.
When nodes coincide, we use the continuous extension in $a_0,\ldots,a_k$.
\end{definition}

Let $\varphi_1,\varphi_2\in\R[t,x]$, and let $a_0,\ldots,a_k$ be real
nodes with $k\ge0$. By the product rule for divided differences
\cite[Corollary~28]{deBoor},
\begin{equation}\label{eq:dd-product}
 (\varphi_1\varphi_2)[a_0,\ldots,a_k]
 =\sum_{j=0}^k \varphi_1[a_0,\ldots,a_j]\,
                    \varphi_2[a_j,\ldots,a_k].
\end{equation}
For pairwise distinct nodes, induction on $b$ using \eqref{eq:dd-product}
for $x^b=x\cdot x^{b-1}$ gives
\begin{equation}\label{eq:dd-monomial}
 (x^b)[a_0,\ldots,a_k]
 =\sum_{\substack{r_0,\ldots,r_k\in\Z_{\ge0}\\r_0+\cdots+r_k=b-k}}
       a_0^{r_0}\cdots a_k^{r_k}.
\end{equation}
The sum in \eqref{eq:dd-monomial} is empty when $b<k$. For $b\ge k$,
every monomial of total degree $b-k$ in $a_0,\ldots,a_k$ occurs exactly
once. Thus the right-hand side of \eqref{eq:dd-monomial} is a symmetric
polynomial in the nodes, and continuity extends \eqref{eq:dd-monomial}
to repeated nodes. Every $\varphi\in\R[t,x]$ is an $\R[t]$-linear
combination of powers of $x$. By linearity of divided differences,
$\varphi[a_0,\ldots,a_k]$ is therefore symmetric in $a_0,\ldots,a_k$.

For each integer $k\ge1$ and $a<b$ with $a\in\R$ and $b\in\R\cup\{+\infty\}$, write
\[
 \Delta_k(a,b)\coloneqq
 \{(v_1,\ldots,v_k)\in(a,b)^k:v_1<\cdots<v_k\}.
\]

\begin{lemma}\label{lem:integral-chains}
Let $\lambda$ be an auxiliary function and $g\in\mathscr F[t]$.
For $\varphi(t,x)\in\R[t,x]$ and $u\in\R$,
\begin{equation}\label{eq:chain-expansion}
\begin{aligned}
 \bigl(\varphi(t,\mathsf X_\lambda)g\bigr)(u)
 &=\varphi(t,\lambda(u))g(u)\\
 &\quad+\sum_{k=1}^{\deg_{x}\varphi}t^k
   \int_{\Delta_k(u,+\infty)}
       \varphi[\lambda(u),\lambda(v_1),\ldots,\lambda(v_k)]g(v_k)
       \,\mathrm{d}v_1\cdots \mathrm{d}v_k .
\end{aligned}
\end{equation}
\end{lemma}

\begin{proof}
Since both sides of \eqref{eq:chain-expansion} are linear in $\varphi$
over $\R[t]$, it suffices to prove the formula for $\varphi(t,x)=x^b$,
where $b\ge0$ is an integer. For $b=0$, both sides equal $g(u)$.
For $b\ge1$, applying $\mathsf X_\lambda$ to $\mathsf X_\lambda^{b-1}g$ gives
\[
 (\mathsf X_\lambda^b g)(u)
 =\lambda(u)(\mathsf X_\lambda^{b-1}g)(u)
  +t\int_u^{+\infty}(\mathsf X_\lambda^{b-1}g)(v)\,\mathrm{d}v.
\]
Choosing multiplication at all $b$ steps keeps the argument at $u$
and gives the unique term without integration, $\lambda(u)^b g(u)$.

Consider a summand of $\mathsf X_\lambda^b g$ with $k\ge1$ integrations.
Write $v_0=u$ and label the integration variables $v_1,\ldots,v_k$ from
the outside in. For $0\le i<k$, let $r_i$ count the multiplications at
$v_i$ before integration in $v_{i+1}$, and let $r_k$ count those at
$v_k$. This choice of integration and multiplication steps gives
\[
 t^k\lambda(u)^{r_0}
 \int_u^{+\infty}\lambda(v_1)^{r_1}\cdots
 \int_{v_{k-1}}^{+\infty}\lambda(v_k)^{r_k}g(v_k)
 \,\mathrm{d}v_k\cdots \mathrm{d}v_1.
\]
The integration region is $\Delta_k(u,+\infty)$. Each tuple of
nonnegative integers with $r_0+\cdots+r_k=b-k$ occurs once, since it
prescribes $r_0$ multiplications, one integration, $r_1$ multiplications,
and so on, ending with $r_k$ multiplications.
By \eqref{eq:dd-monomial}, summing $\prod_{i=0}^k\lambda(v_i)^{r_i}$
over these tuples gives $(x^b)[\lambda(v_0),\ldots,\lambda(v_k)]$.
Adding the resulting integrals for $1\le k\le b$ to $\lambda(u)^b g(u)$
gives \eqref{eq:chain-expansion} for $x^b$.
\end{proof}

To identify the kernel in \eqref{eq:chain-expansion}, write $a=u$ and
take $v=v_k$ as the outer variable. For $2\le k\le\deg_x\varphi$, the
$k$th summand becomes
\[
 \int_a^{+\infty}
 \underbrace{\left(
   t^k\!\int_{\Delta_{k-1}(a,v)}
   \varphi[\lambda(a),\lambda(v_1),\ldots,\lambda(v_{k-1}),\lambda(v)]
   \,\mathrm{d}v_1\cdots \mathrm{d}v_{k-1}
 \right)}_{\text{the $k$th summand of the kernel}}
 g(v)\,\mathrm{d}v .
\]

Choose $M>a$ so that $g(v)=0$ for $v>M$. By local boundedness,
each coefficient of $g(v)$ in $t$ is bounded on $[a,M]$.
Taking the coefficient of $t^\ell$ commutes with taking a divided
difference in $x$. Apply \eqref{eq:dd-monomial} to each power of $x$
in $[t^\ell]\varphi(t,x)$. The resulting polynomial in
$\lambda(a),\lambda(v_1),\ldots,\lambda(v_{k-1}),\lambda(v)$
is a finite sum of scalar multiples of monomials. Each node lies
in $[0,1]$, so each monomial has absolute value at most one.
The sum of the absolute values of the scalar coefficients therefore
gives a bound independent of $v_1,\ldots,v_{k-1},v$.
The integrand vanishes for $v>M$, and each of its $t$-coefficients is
bounded on the finite-volume region $a<v_1<\cdots<v_{k-1}<v\le M$.
Fubini's theorem thus applies to each coefficient.

Fix $v>a$ and put $m=k-1$. Sorting $v_1,\ldots,v_m$ gives
\[
 (a,v)^m\setminus\bigcup_{i<j}\{v_i=v_j\}
 =\bigsqcup_{(i_1,\ldots,i_m)}
   \{a<v_{i_1}<\cdots<v_{i_m}<v\},
\]
where $(i_1,\ldots,i_m)$ runs through the permutations of $(1,\ldots,m)$.
For example, when $m=2$, swapping $v_1$ and $v_2$ exchanges the two
triangles of the square. Since permutations of $v_1,\ldots,v_m$ preserve
Lebesgue measure and the integrand
$\varphi[\lambda(a),\lambda(v_1),\ldots,\lambda(v_m),\lambda(v)]$
is symmetric in its interior nodes, the integrals over the $m!$ ordered
regions are equal. The diagonals $v_i=v_j$ with $i<j$ have measure zero, so the integral
over $a<v_1<\cdots<v_m<v$ is $1/m!$ times the integral over $(a,v)^m$.

In the integral over $(a,v)^m$, the divided difference depends on
each $v_i$ only through $\lambda(v_i)$. Let $\rho_{a,v}$ be the
pushforward of Lebesgue measure on $(a,v)$ under $\lambda$.
Applying \eqref{eq:occupation-values} of Appendix~\ref{app:measure-conventions}
to one variable $v_i$ at a time, keeping the others fixed, gives
\begin{equation}\label{eq:interior-pushforward}
\begin{aligned}
 &\frac{t^{m+1}}{m!}\int_{(a,v)^m}
 \varphi[\lambda(a),\lambda(v_1),\ldots,\lambda(v_m),\lambda(v)]\,
 \mathrm{d}v_1\cdots \mathrm{d}v_m\\
 &\quad=\frac{t^{m+1}}{m!}\int_{[0,1]^m}
 \varphi[\lambda(a),z_1,\ldots,z_m,\lambda(v)]\,
 \mathrm{d}\rho_{a,v}^{\otimes m}(z_1,\ldots,z_m).
\end{aligned}
\end{equation}
For $m=0$, the right-hand side of \eqref{eq:interior-pushforward}
gives the first kernel summand $t\varphi[\lambda(a),\lambda(v)]$.

For fixed $z_1,\ldots,z_m$ and $\lambda(a)\ne\lambda(v)$, the
recursion in Definition~\ref{def:divided-differences} and symmetry
in the nodes give
\[
 \varphi[\lambda(a),z_1,\ldots,z_m,\lambda(v)]
 =\frac{\varphi[\lambda(v),z_1,\ldots,z_m]
        -\varphi[\lambda(a),z_1,\ldots,z_m]}
       {\lambda(v)-\lambda(a)}.
\]
For $\lambda(a)=\lambda(v)$, continuous extension replaces the
quotient by the derivative of $x\mapsto\varphi[x,z_1,\ldots,z_m]$
at $x=\lambda(a)$. For all endpoint values, linearity allows us
to integrate $\varphi[x,z_1,\ldots,z_m]$ over the interior nodes
first and then take the first divided difference of the resulting
polynomial in $x$ at $\lambda(a)$ and $\lambda(v)$.

\begin{definition}\label{def:occupation-transport}
For a finite Borel measure $\rho$ on $[0,1]$, define the $\R[t]$-linear
operator $\mathsf T_\rho:\R[t,x]\longrightarrow\R[t,x]$ by setting,
for each $\varphi\in\R[t,x]$,
\begin{equation}\label{eq:occupation-transport}
 (\mathsf T_\rho \varphi)(t,x)
 \coloneqq\sum_{j=0}^{\deg_x\varphi}\frac{t^j}{j!}
   \int_{[0,1]^j}\varphi[x,z_1,\ldots,z_j]\,
                    \mathrm{d}\rho^{\otimes j}(z).
\end{equation}
The term with $j=0$ is $\varphi(t,x)$.
For real nodes $x_0,x_1$, write
$(\mathsf T_\rho\varphi)[x_0,x_1]$ for the first divided difference
of the polynomial $(\mathsf T_\rho\varphi)(t,x)$ in $x$ at $x_0,x_1$.
The letter $\mathsf T$ stands for \emph{transport} under interval
subdivision (Lemma~\ref{lem:occupation-semigroup}).
\end{definition}

With $\rho=\rho_{a,v}$, the weight $t^{j+1}/j!$ in
\eqref{eq:interior-pushforward} for $m=j$ is $t$ times the weight
$t^j/j!$ in \eqref{eq:occupation-transport}.
The summand with $j=\deg_x\varphi$ is constant in $x$ and has zero
first divided difference. Thus the integral part of
$\bigl(\varphi(t,\mathsf X_\lambda)g\bigr)(a)$ is
\[
 t\int_a^{+\infty}
   (\mathsf T_{\rho_{a,v}}\varphi)[\lambda(a),\lambda(v)]\,g(v)\,\mathrm{d}v.
\]

\begin{lemma}\label{lem:occupation-semigroup}
For finite Borel measures $\rho,\widetilde\rho$ on $[0,1]$,
\begin{equation}\label{eq:occupation-semigroup}
 \mathsf T_{\rho+\widetilde\rho}=\mathsf T_\rho\mathsf T_{\widetilde\rho}.
\end{equation}
\end{lemma}

\begin{proof}
Fix $\varphi\in\R[t,x]$. By \eqref{eq:occupation-transport}, the
$j$th summand of $(\mathsf T_{\widetilde\rho}\varphi)(t,\xi)$ is
\[
 \frac{t^j}{j!}\int_{[0,1]^j}
 \varphi[\xi,w_1,\ldots,w_j]\,
 \mathrm{d}\widetilde\rho^{\otimes j}(w_1,\ldots,w_j).
\]
Integration over $w_1,\ldots,w_j$ leaves a polynomial in $\xi$.
In the $i$th summand for the outer operator $\mathsf T_\rho$, we take
the divided difference of this polynomial at $x,z_1,\ldots,z_i$.
By linearity, we may take the divided difference before integrating in
$w_1,\ldots,w_j$.
The composition rule \cite[equation~(18)]{deBoor} and symmetry in the nodes give
\begin{equation}\label{eq:dd-composition}
 \bigl(\xi\mapsto\varphi[\xi,w_1,\ldots,w_j]\bigr)[x,z_1,\ldots,z_i]
 =\varphi[x,z_1,\ldots,z_i,w_1,\ldots,w_j].
\end{equation}
Substituting \eqref{eq:dd-composition} into the expansion of
$\mathsf T_\rho(\mathsf T_{\widetilde\rho}\varphi)$ gives
\[
 (\mathsf T_\rho\mathsf T_{\widetilde\rho}\varphi)(t,x)
 =\sum_{i,j\ge0}\frac{t^{i+j}}{i!j!}
   \int_{[0,1]^i\times[0,1]^j}\!\varphi[x,z_1,\ldots,z_i,w_1,\ldots,w_j]\,
           \mathrm{d}\rho^{\otimes i}(z)\,\mathrm{d}\widetilde\rho^{\otimes j}(w).
\]
Group the terms by $n=i+j$, omitting $n>\deg_x\varphi$ since the
divided difference of order $n$ then vanishes.
In $(\rho+\widetilde\rho)^{\otimes n}$ there are $\binom ni$ choices
of the $i$ coordinates carrying $\rho$, with $\widetilde\rho$ on the
remaining coordinates. Moving the chosen coordinates to the first $i$
positions preserves the integrand by symmetry of the divided difference.
Thus all $\binom ni$ choices give the same integral against
$\rho^{\otimes i}\otimes\widetilde\rho^{\otimes(n-i)}$.
Since $\binom ni t^n/n!=t^{i+j}/(i!j!)$ for $i+j=n$, the sum over
all pairs with total $n$ is the $n$th summand of
$(\mathsf T_{\rho+\widetilde\rho}\varphi)(t,x)$ in \eqref{eq:occupation-transport}.
\end{proof}

For $a<w<v$, splitting $(a,v)$ at $w$ gives $\rho_{a,v}=\rho_{a,w}+\rho_{w,v}$
because the single position $w$ has Lebesgue measure zero.
Lemma~\ref{lem:occupation-semigroup} then gives
$\mathsf T_{\rho_{a,v}}=\mathsf T_{\rho_{a,w}}\mathsf T_{\rho_{w,v}}$.

For fixed $u$, changing $K_{\boldsymbol{\varphi},\lambda}(u,v)$ on a set of
Lebesgue measure zero in $v$ leaves the integral against every test function
$g$ unchanged. Lemma~\ref{lem:boundary-limits} uses the pointwise values
of the kernel at $(u,v)=(0,h+s)$ for integers $h\ge0$ and at $(0,h-s)$
for integers $h\ge1$, with $0<s<1$.

\begin{definition}\label{def:canonical-kernel}
Let $\boldsymbol{\varphi}=(\varphi_j)_{j\ge0}$ be a finite normal form
and $\lambda$ an auxiliary function.
Fix the evaluation point $u\in\R$, write $v$ for the integration
variable, and call $\Delta=v-u$ the \emph{jump}.
For $\Delta>0$ with $\Delta\notin\Z$, define the \emph{integral kernel} by
\begin{equation}\label{eq:2-14}
 K_{\boldsymbol{\varphi},\lambda}(u,u+\Delta)
 \coloneqq t\sum_{0\le j<\Delta}
   (\mathsf T_{\rho_{u+j,u+\Delta}}\varphi_j)[\lambda(u),\lambda(u+\Delta)].
\end{equation}
Set $K_{\boldsymbol{\varphi},\lambda}(u,u+\Delta)=0$ for $\Delta\in\Z$ or $\Delta\le0$.
For $g\in\mathscr F[t]$ and $h\in\Z_{\ge0}$,
$(\mathsf R_{\lambda,\boldsymbol\varphi}g)(u)$ includes the separate
value $g(u+h)$ with \emph{evaluation weight} $\varphi_h(t,\lambda(u))$.
For $a\in\R$, let $\operatorname{Dirac}(a)$ denote the unit point mass at $a$,
so $\operatorname{Dirac}(a)(B)=\mathbf1_B(a)$ for every Borel set $B\subseteq\R$.
The atom at $v=u+h$ is $\varphi_h(t,\lambda(u))\operatorname{Dirac}(u+h)$.
\end{definition}

The shift $\mathsf S^j$ starts the expansion \eqref{eq:chain-expansion}
of the $j$th component at $u+j$. Since $\lambda(u+j)=\lambda(u)$,
the term without integration is $\varphi_j(t,\lambda(u))g(u+j)$.
For positive noninteger $\Delta$, an integration chain from $u+j$ to
$v=u+\Delta$ requires $u+j<v$, or equivalently $j<\Delta$.
Summing $t(\mathsf T_{\rho_{u+j,v}}\varphi_j)[\lambda(u),\lambda(v)]$
over $0\le j<\Delta$ gives \eqref{eq:2-14}.
The factor $t$ in \eqref{eq:2-14} and the zero value at integer or
nonpositive displacements give $K_{\boldsymbol{\varphi},\lambda}(u,v)\in t\R[t]$
for all $u,v\in\R$.

For $L>0$ and each integer $k\ge0$, we claim that
$[t^k]K_{\boldsymbol\varphi,\lambda}(u,v)$ is jointly Borel measurable
on $0<v-u\le L$ and bounded there by a constant depending only on
$\boldsymbol\varphi,k,L$. For a chain from $u+j$ to $v$ with $m\ge1$
interior positions, put $w_i=v_i-u-j$. In the fixed box $[0,L]^m$,
the indicator of $0<w_1<\cdots<w_m<v-u-j$ is Borel in
$(u,v,w_1,\ldots,w_m)$. By \eqref{eq:dd-monomial}, each coefficient in $t$ of
$\varphi_j[\lambda(u),\lambda(u+j+w_1),\ldots,\lambda(u+j+w_m),\lambda(v)]$
is a polynomial in the $m+2$ displayed values of the Borel function $\lambda$.
Integration over the box with the indicator therefore gives a Borel
function of $(u,v)$. The polynomial is bounded on $[0,1]^{m+2}$,
and the ordered region has volume at most $L^m/m!$.
For $m=0$ only the endpoint nodes in $\varphi_j[\lambda(u),\lambda(v)]$
remain. For noninteger $v-u$, formula \eqref{eq:2-14} adds the chain
integrals for the nonzero components $\varphi_j$ with $0\le j<v-u$,
each with chain lengths $1,\ldots,\deg_x\varphi_j$.
On the remaining lines $v-u=n$ for integers $1\le n\le L$,
Definition~\ref{def:canonical-kernel} gives $[t^k]K_{\boldsymbol\varphi,\lambda}(u,v)=0$.

Changing finitely many values of $\lambda$ per period preserves the
occupation measures $\rho_{u+j,v}$ in \eqref{eq:2-14}, but may change
either node $\lambda(u)$ or $\lambda(v)$ of the first divided difference.
For example, take $\boldsymbol{\varphi}=\boldsymbol\alpha_2^2$.
For $u<v<u+1$, only the component $x^2$ at $j=0$ enters
\eqref{eq:2-14}. Since $[t^0]\mathsf T_{\rho_{u,v}}(x^2)=x^2$
by \eqref{eq:occupation-transport},
\[
 [t]K_{\boldsymbol\alpha_2^2,\lambda}(u,v)
 =(x^2)[\lambda(u),\lambda(v)]=\lambda(u)+\lambda(v).
\]
Let $\widetilde\lambda$ be an auxiliary function agreeing with $\lambda$
outside $u+\Z$. For every $v\in(u,u+1)$ we have
$\widetilde\lambda(v)=\lambda(v)$, so only the first node can change and
\[
 [t]K_{\boldsymbol\alpha_2^2,\widetilde\lambda}(u,v)
 -[t]K_{\boldsymbol\alpha_2^2,\lambda}(u,v)
 =\widetilde\lambda(u)-\lambda(u).
\]

Fix finite normal forms $\boldsymbol{\varphi}=(\varphi_a)_{a\ge0}$ and
$\boldsymbol{\psi}=(\psi_b)_{b\ge0}$, an auxiliary function $\lambda$,
and an initial position $u\in\R$. Write
$K_1\coloneqq K_{\boldsymbol{\varphi},\lambda}$ and
$K_2\coloneqq K_{\boldsymbol{\psi},\lambda}$.

\begin{samepage}
\begin{lemma}\label{lem:mixed-kernel-composition}
The operators satisfy
$\mathsf R_{\lambda,\boldsymbol\varphi\boldsymbol\psi}
=\mathsf R_{\lambda,\boldsymbol\varphi}\mathsf R_{\lambda,\boldsymbol\psi}$.
For $g\in\mathscr F[t]$ and an integer $h\ge0$, the evaluation term at
$u+h$ in $(\mathsf R_{\lambda,\boldsymbol\varphi\boldsymbol\psi}g)(u)$ is
$(\boldsymbol{\varphi}\boldsymbol{\psi})_h(t,\lambda(u))g(u+h)$.
For $v>u$ with $v-u\notin\Z$, the integral kernel of
$\mathsf R_{\lambda,\boldsymbol\varphi\boldsymbol\psi}$ is
\begin{equation}\label{eq:2-16}
 \begin{aligned}
 K_{\boldsymbol{\varphi}\boldsymbol{\psi},\lambda}(u,v)
 &=\sum_{a\ge0}\varphi_a(t,\lambda(u))K_2(u+a,v)
   +\sum_{b\ge0}K_1(u,v-b)\psi_b(t,\lambda(v))\\
 &\quad+\int_u^vK_1(u,w)K_2(w,v)\,\mathrm{d}w.
 \end{aligned}
\end{equation}
\end{lemma}
\end{samepage}

\begin{proof}
Expand both factors in
$\mathsf R_{\lambda,\boldsymbol\varphi}\mathsf R_{\lambda,\boldsymbol\psi}$
by \eqref{eq:normal-form-operator}. Commuting $\mathsf S^b$ past
$\varphi_a(t,\mathsf X_\lambda)$ and grouping the products by $a+b$ gives
\[
 \mathsf R_{\lambda,\boldsymbol\varphi}\mathsf R_{\lambda,\boldsymbol\psi}
 =\sum_{a,b\ge0}\mathsf S^{a+b}(\varphi_a\psi_b)(t,\mathsf X_\lambda)
 =\mathsf R_{\lambda,\boldsymbol\varphi\boldsymbol\psi}.
\]
The coefficient of the separate evaluation $g(u+h)$ is
$\sum_{a+b=h}\varphi_a(t,\lambda(u))\psi_b(t,\lambda(u+a))
=(\boldsymbol\varphi\boldsymbol\psi)_h(t,\lambda(u))$
by periodicity of $\lambda$.

Fix $v>u$ with $v-u\notin\Z$. For integer shifts $a,b\ge0$ with $a+b<v-u$,
consider a chain $u+a+b=z_0<\cdots<z_k=v$ with $k\ge1$ in the expansion
\eqref{eq:chain-expansion} for $\varphi_a\psi_b$.
In the product rule \eqref{eq:dd-product} for
$(\varphi_a\psi_b)[\lambda(z_0),\ldots,\lambda(z_k)]$,
the summand indexed by $i\in\{0,\ldots,k\}$ is
$\varphi_a[\lambda(z_0),\ldots,\lambda(z_i)]\,
 \psi_b[\lambda(z_i),\ldots,\lambda(z_k)]$.
Put $w=z_i-b$. The positions in the two integration chains are
\[
 \underbrace{(z_0-b,\ldots,z_i-b)}_{\varphi_a},
 \qquad
 \underbrace{(z_i,\ldots,z_k)}_{\psi_b}.
\]
Translation by $-b$ preserves Lebesgue measure and, by periodicity,
all node values $\lambda(z_j-b)=\lambda(z_j)$ for $0\le j\le i$.
Adding $b$ to the positions in the $\varphi_a$ chain and joining the
two chains at $z_i$ recovers $(z_0,\ldots,z_k)$ with $z_i$ marked.

For $i=0$ or $i=k$, only one chain has positive length. By \eqref{eq:chain-expansion}
and \eqref{eq:2-14}, we integrate the divided difference on that chain
over its $k-1$ ordered interior positions and sum with weight $t^k$
over $k\ge1$ and the shift of the integrated factor.
For $i=0$, fix $0\le a<v-u$. The $\psi_b$ chain has interior positions in
$(u+a+b,v)$ with $0\le b<v-u-a$, giving
$\varphi_a(t,\lambda(u))K_2(u+a,v)$ since $\lambda(u+a)=\lambda(u)$.
For $i=k$, fix $0\le b<v-u$. The $\varphi_a$ chain has interior positions in
$(u+a,v-b)$ with $0\le a<v-u-b$, giving $K_1(u,v-b)\psi_b(t,\lambda(v))$.
Summing over the fixed shifts $a$ and $b$, respectively, gives the first
two sums in \eqref{eq:2-16}.

For $0<i<k$, fix $w$ with $u+a<w<v-b$. The positions $z_j-b$ for
$0<j<i$ range in increasing order between $u+a$ and $w$,
independently of the positions $z_j$ for $i<j<k$ between $w+b$ and $v$.
The two chain lengths $i$ and $k-i$ give the factors $t^i$ and $t^{k-i}$.
By Fubini's theorem, we may integrate over the two sets of interior
positions before integrating over $w=z_i-b$. If $w-u\notin\Z$ and
$v-w\notin\Z$, summing the weighted integrals over both positive
lengths and over $0\le a<w-u$, $0\le b<v-w$ gives
$K_1(u,w)K_2(w,v)$ by \eqref{eq:2-14}.
For fixed $u,v$, only finitely many $w\in(u,v)$ have $w-u\in\Z$ or
$v-w\in\Z$, so integration over $w$ gives the last term of \eqref{eq:2-16}.
Thus \eqref{eq:2-16} holds for every $v>u$ with $v-u\notin\Z$ and uses
the prescribed values $\lambda(u)$ and $\lambda(v)$.
\end{proof}

\subsection{Closed truncation and kernel limits}\label{subsec:closed-truncation}

Fix a finite normal form $\boldsymbol{\varphi}$ of $\beta(t,z)$ as in \eqref{eq:2-1}
and an integer $h\ge0$. For each positive integer $t$, define
the truncated generating series
\[
\begin{aligned}
 \beta^{\le h}(t,z)&\coloneqq\sum_{a=0}^{ht}([z^a]\beta(t,z))z^a,\\
 \beta^{\ge h}(t,z)&\coloneqq\sum_{a=ht}^{+\infty}([z^a]\beta(t,z))z^a.
\end{aligned}
\]
The boundary term $\bigl([z^{ht}]\beta(t,z)\bigr)z^{ht}$ is retained
in both $\beta^{\le h}(t,z)$ and $\beta^{\ge h}(t,z)$.

We seek correction polynomials $\delta_h^+$ for the upper bound $a\le ht$
and $\delta_h^-$ for the lower bound $a\ge ht$, both in $\R[t,x]$.
With $x=(1-z)^{-1}$, we require for every positive integer $t$ that
\[
\begin{aligned}
 \beta^{\le h}(t,z)
   &=\sum_{j\le h}z^{jt}\varphi_j(t,x)+z^{ht}\delta_h^+(t,x),\\
 \beta^{\ge h}(t,z)
   &=\sum_{j\ge h}z^{jt}\varphi_j(t,x)+z^{ht}\delta_h^-(t,x).
\end{aligned}
\]
Keeping only the components with $j\ge h$ omits any terms of exponent
$a\ge ht$ from components with $j<h$, so $z^{ht}\delta_h^-(t,x)$ must
restore those terms in $\beta^{\ge h}(t,z)$.

Let $\eta$ be Lebesgue measure on $[0,1]$, and put
$\mathsf T_\ell\coloneqq\mathsf T_{\ell\eta}$ for $\ell\ge0$.
By periodicity and Definition~\ref{def:auxiliary-function}(ii),
$\rho_{a,a+1}=\eta$ for every $a\in\R$.
Define the boundary corrections by
\begin{equation}\label{eq:2-10}
\begin{aligned}
 \delta_h^+(t,x)&\coloneqq(1-x)\sum_{j\le h}(\mathsf T_{h-j}\varphi_j)[x,1],\\
 \delta_h^-(t,x)&\coloneqq x\sum_{j<h}(\mathsf T_{h-j}\varphi_j)[0,x].
\end{aligned}
\end{equation}
The parameter $h-j=(ht-jt)/t$ is the distance from $jt$ to the
cutoff exponent $ht$ in units of $t$.
The definitions give $\delta_h^+(t,1)=\delta_h^-(t,0)=0$
and $\delta_0^-=0$.

Define finite polynomial sequences
$\boldsymbol{\varphi}^{\le h}=(\varphi_j^{\le h}(t,x))_{j\ge0}$ and
$\boldsymbol{\varphi}^{\ge h}=(\varphi_j^{\ge h}(t,x))_{j\ge0}$ by
\[
\begin{aligned}
 \varphi_j^{\le h}(t,x)&\coloneqq
 \begin{cases}
  \varphi_j(t,x),&j<h,\\
  \varphi_h(t,x)+\delta_h^+(t,x),&j=h,\\
  0,&j>h,
 \end{cases}\\[6pt]
 \varphi_j^{\ge h}(t,x)&\coloneqq
 \begin{cases}
  0,&j<h,\\
  \varphi_h(t,x)+\delta_h^-(t,x),&j=h,\\
  \varphi_j(t,x),&j>h.
 \end{cases}
\end{aligned}
\]

To prove that the corrections in \eqref{eq:2-10} give the closed truncations
$\beta^{\le h}$ and $\beta^{\ge h}$, we first compute $\mathsf T_{h-j}$
on the monomials of $\varphi_j$.

\begin{lemma}\label{lem:uniform-traces}
For $\ell\ge0$ and an integer $b\ge0$,
\begin{equation}\label{eq:uniform-transport-monomial}
 \mathsf T_\ell(x^b)
 =\sum_{j=0}^b\binom{\ell t+j-1}{j}x^{b-j}.
\end{equation}
\end{lemma}

\begin{proof}
Fix an integer $0\le j\le b$. In the summand of \eqref{eq:occupation-transport}
with $k$ integration variables, \eqref{eq:dd-monomial} expands
$(x^b)[x,z_1,\ldots,z_k]$ as the sum of all monomials
$x^{r_0}z_1^{r_1}\cdots z_k^{r_k}$ of total degree $b-k$.
To extract $x^{b-j}$, take $r_0=b-j$, so $r_1+\cdots+r_k=j-k$.
The positive integers $l_i=r_i+1$ then sum to $j$, and integrating
$z_i^{r_i}$ against $\ell\,\mathrm{d}z_i$ on $[0,1]$ gives $\ell/l_i$.
Hence
\[
 [x^{b-j}]\mathsf T_\ell(x^b)
 =\sum_{k=0}^j\frac{(\ell t)^k}{k!}
   \sum_{\substack{l_1+\cdots+l_k=j\\l_i\ge1}}
        \frac1{l_1\cdots l_k}.
\]

Arrange the labels $1,\ldots,j$ in a row, cut the row into successive
blocks of prescribed lengths $l_1,\ldots,l_k$, and read each block as
a cycle. Each cycle of
length $l_i$ has $l_i$ choices of starting label, so the $j!$ rows give
$j!/(l_1\cdots l_k)$ distinct ordered lists of cycles.
After summing over all ordered length tuples, we divide by $k!$ because
each permutation with $k$ cycles occurs once for each of the $k!$
orderings of its cycles. For $j=0$, the empty permutation contributes $1$. Let $c(\tau)$ denote the number of cycles of a permutation
$\tau\in\mathfrak S_j$. By the cycle-enumerator identity
\cite[Proposition~1.3.7]{StanleyEC1},
\[
 [x^{b-j}]\mathsf T_\ell(x^b)
 =\frac{1}{j!}\sum_{\tau\in\mathfrak S_j}(\ell t)^{c(\tau)}
 =\binom{\ell t+j-1}{j}.\qedhere
\]
\end{proof}

For $\ell\ge0$ and an integer $b\ge1$, taking the first divided difference
at $x,1$ in \eqref{eq:uniform-transport-monomial} replaces each $x^{b-j}$
by $\sum_{i=0}^{b-j-1}x^i$. To obtain the coefficient of $x^i$, we therefore
sum the binomial coefficients over $0\le j\le b-1-i$.
\begin{equation}\label{eq:uniform-wall-polynomial}
 (\mathsf T_\ell x^b)[x,1]
 =\sum_{i=0}^{b-1}\left(\sum_{j=0}^{b-1-i}\binom{\ell t+j-1}{j}\right)x^i
 =\sum_{i=0}^{b-1}\binom{\ell t+b-1-i}{b-1-i}x^i.
\end{equation}
Setting $x=0$ in \eqref{eq:uniform-wall-polynomial} gives
\begin{equation}\label{eq:uniform-wall-count}
 (\mathsf T_\ell x^b)[0,1]=\binom{\ell t+b-1}{b-1}.
\end{equation}
For $b=0$, both divided differences vanish because $\mathsf T_\ell 1=1$.

\begin{samepage}
\begin{lemma}\label{lem:closed-clipping}
Let $\boldsymbol{\varphi}$ be a finite normal form of $\beta(t,z)$,
and let $h\ge0$ be an integer.
Then $\boldsymbol{\varphi}^{\le h}$ and
$\boldsymbol{\varphi}^{\ge h}$ are finite normal forms of
$\beta^{\le h}$ and $\beta^{\ge h}$, respectively.
For every auxiliary function $\lambda$, every $u\in\R$,
and every noninteger $\Delta>0$,
\[
\begin{aligned}
 K_{\boldsymbol{\varphi}^{\le h},\lambda}(u,u+\Delta)&=\mathbf1_{\Delta<h}K_{\boldsymbol{\varphi},\lambda}(u,u+\Delta),\\
 K_{\boldsymbol{\varphi}^{\ge h},\lambda}(u,u+\Delta)&=\mathbf1_{\Delta>h}K_{\boldsymbol{\varphi},\lambda}(u,u+\Delta).
\end{aligned}
\]
\end{lemma}
\end{samepage}

\begin{proof}
Define the polynomial $\chi_h\in\R[t,x]$ by
\begin{equation}\label{eq:wall-polynomial}
 \chi_h\coloneqq\sum_{j=0}^{h-1}\mathsf T_{h-j}\varphi_j,
\end{equation}
where the empty sum for $h=0$ gives $\chi_0=0$.
By \eqref{eq:2-10} and $\mathsf T_0\varphi_h=\varphi_h$,
$\delta_h^-=x\chi_h[0,x]$ and
$\delta_h^+=(1-x)(\chi_h+\varphi_h)[x,1]$.
Expanding the two divided differences gives
\begin{equation}\label{eq:clipping-local-identities}
\begin{aligned}
 \delta_h^-(t,x)&=\chi_h(t,x)-\chi_h(t,0),\\
 \varphi_h(t,x)+\delta_h^+(t,x)
   &=\varphi_h(t,1)+\chi_h(t,1)-\chi_h(t,x).
\end{aligned}
\end{equation}

Fix a positive integer $t$ and put $x=(1-z)^{-1}$. For $j<h$ and
$b\ge1$, set $k=(h-j)t$. The Chaundy--Bullard identity
\cite[equation~(1.6)]{KS} with parameters $k,b-1$ and argument $z$,
together with \eqref{eq:uniform-wall-polynomial} and $1-x=-zx$, gives
\[
 \sum_{a=0}^k\binom{a+b-1}{b-1}z^a
 =x^b+z^k\bigl((\mathsf T_{h-j}x^b)(t,1)
                     -(\mathsf T_{h-j}x^b)(t,x)\bigr).
\]
Write $\varphi_j(t,x)=\sum_{b\ge0}\alpha_{jb}(t)x^b$.
Multiply the monomial identity by $z^{jt}\alpha_{jb}(t)$ and sum over
$b\ge1$. For $b=0$, the exponent $jt$ is below $ht$, and the correction
vanishes because $\mathsf T_{h-j}$ fixes constants in $x$. Since $jt+k=ht$,
the upper truncation of $z^{jt}\varphi_j(t,x)$ for $j<h$ is
\[
 z^{jt}\varphi_j(t,x)
 +z^{ht}\bigl((\mathsf T_{h-j}\varphi_j)(t,1)
             -(\mathsf T_{h-j}\varphi_j)(t,x)\bigr).
\]
For $j=h$, only $z^{ht}\varphi_h(t,1)$ remains, since
$[z^0]\varphi_h(t,(1-z)^{-1})=\varphi_h(t,1)$.
For $j>h$, every exponent is at least $jt>ht$.
Summing over $j$ and using \eqref{eq:clipping-local-identities}
gives the normal form $\boldsymbol{\varphi}^{\le h}$ of $\beta^{\le h}$.

Subtracting $\beta^{\le h}$ from $\beta$ leaves only the terms of
exponent greater than $ht$. Restoring the term of exponent $ht$ gives
$\beta^{\ge h}=\beta-\beta^{\le h}+\zeta_h(t)z^{ht}$.
By \eqref{eq:uniform-wall-count}, the sum over $j<h$ in the coefficient
formula \eqref{eq:2-3} is $\chi_h[0,1]=\chi_h(t,1)-\chi_h(t,0)$.
Adding the two identities in \eqref{eq:clipping-local-identities} now gives
$\varphi_h+\delta_h^++\delta_h^-=\zeta_h(t)$. Thus for every $j\ge0$,
\[
 \varphi_j^{\ge h}
 =\varphi_j-\varphi_j^{\le h}+\mathbf1_{\{j=h\}}\zeta_h(t).
\]
Hence $\boldsymbol{\varphi}^{\ge h}$ represents $\beta^{\ge h}$.

For a noninteger $\Delta$ with $0<\Delta<h$, the upper truncation retains every
component with $j<\Delta$ in the kernel formula \eqref{eq:2-14}.
For a noninteger $\Delta>h$, put $w=u+h$ and $v=u+\Delta$.
The interval $(u+j,w)$ contains $h-j$ complete periods for $j<h$,
so $\rho_{u+j,v}=(h-j)\eta+\rho_{w,v}$. By
\eqref{eq:occupation-semigroup} and \eqref{eq:2-14},
\[
 \frac{K_{\boldsymbol{\varphi}^{\le h},\lambda}(u,v)}t
 =\bigl(\mathsf T_{\rho_{w,v}}(\chi_h+\varphi_h+\delta_h^+)\bigr)
       [\lambda(u),\lambda(v)].
\]
By \eqref{eq:clipping-local-identities},
$\chi_h(t,x)+\varphi_h(t,x)+\delta_h^+(t,x)=\varphi_h(t,1)+\chi_h(t,1)$.
The right-hand side is constant in $x$, so $\mathsf T_{\rho_{w,v}}$
fixes it and its first divided difference is zero for all endpoint
values $\lambda(u),\lambda(v)$. Hence
$K_{\boldsymbol{\varphi}^{\le h},\lambda}(u,v)=0$ for $\Delta>h$.

The boundary term at $j=h$ adds no integral kernel because
$\mathsf T_\rho\zeta_h(t)=\zeta_h(t)$ is constant in $x$ and has zero
first divided difference. Applying \eqref{eq:2-14} to the component
identity for $\varphi_j^{\ge h}$ gives
$K_{\boldsymbol{\varphi}^{\ge h},\lambda}
=K_{\boldsymbol{\varphi},\lambda}-K_{\boldsymbol{\varphi}^{\le h},\lambda}$
at every noninteger $\Delta>0$.
\end{proof}

For every positive integer $t$, the lower truncation retains the coefficient
$[z^{ht}]\beta(t,z)$. In the normal form of $\beta^{\ge h}$, the components
with $j<h$ vanish, and each summand with $j>h$ starts at degree at least
$jt>ht$ in $z$. Thus $[z^{ht}]\beta(t,z)$ is the constant term in $z$ of
$\varphi_h^{\ge h}(t,(1-z)^{-1})$. At $z=0$ the second argument is $1$,
so uniqueness in Lemma~\ref{lem:coefficient-polynomial} gives
\begin{equation}\label{eq:boundary-coefficient}
 \zeta_h(t)=\varphi_h^{\ge h}(t,1)=\varphi_h(t,1)+\delta_h^-(t,1).
\end{equation}

\begin{example}
\label{ex:running-clipping}
For $I=\{1,2\}$ in Example~\ref{ex:running-polypositroid}, the inequality
$y_1+y_2\le2$ cuts the square $[0,2]^2$ down to the triangle $\mathcal P_I$.
Fix a positive integer $t$ and put $x=(1-z)^{-1}$. For each integer
$0\le a\le 2t$, the lattice points of $t\mathcal P_I$ with coordinate sum
$a$ are $(i,a-i)$ for $i=0,\ldots,a$. Hence $\beta_I(t,z)$ agrees with
$x^2=\sum_{a\ge0}(a+1)z^a$ through degree $2t$. No point of
$t\mathcal P_I$ has a larger coordinate sum. Writing $a=2t+1+i$ in the tail
and using $1-x=-zx$ gives
\[
 \beta_I(t,z)
 =x^2-z^{2t+1}\bigl((2t+1)x+x^2\bigr)
 =x^2+z^{2t}(1-x)(2t+1+x).
\]
The coefficient $[z^{2t}]\beta_I(t,z)=2t+1$ counts the lattice points on
the edge $y_1+y_2=2t$ and comes entirely from $x^2$.

The generating function $\beta_{[0,2]}^2$ of the square has normal form
$\boldsymbol\alpha_2^2$, with components
$x^2$, $2x(1-x)$, and $(1-x)^2$ at indices $0,2,4$. The upper truncation
of Lemma~\ref{lem:closed-clipping} at $h=2$ removes the component at index
$4$ and replaces $2x(1-x)$ at index $2$ by $(1-x)(2t+1+x)$.
The amount added to the old component is therefore
\[
 \delta_2^+(t,x)=(1-x)(2t+1+x)-2x(1-x)=(1-x)(2t+1-x).
\]
For an auxiliary function $\lambda$, a test function $g\in\mathscr F[t]$,
and $u\in\R$, this correction adds the evaluation term
\[
 \delta_2^+(t,\lambda(u))g(u+2)
 =(1-\lambda(u))(2t+1-\lambda(u))g(u+2).
\]
\end{example}

To study the signs of the coefficients in $t$ of boundary corrections such
as $\delta_2^+(t,x)$ in Example~\ref{ex:running-clipping}, we express
$\delta_h^\pm(t,x)$ in terms of one-sided limits of the kernel quotient
$K_{\boldsymbol{\varphi},\lambda}/t$ before truncation.

\begin{lemma}\label{lem:boundary-limits}
Fix $x\in[0,1]$. There exist two families of auxiliary functions
$(\lambda_s^+)_{0<s<1}$ and $(\lambda_s^-)_{0<s<1}$ such that,
for every finite normal form $\boldsymbol{\varphi}$, the upper correction
at each integer $h\ge0$ satisfies
\begin{equation}\label{eq:boundary-upper-limit}
 \delta_h^+(t,x)
 =(1-x)\lim_{s\downarrow0}
       \frac{K_{\boldsymbol{\varphi},\lambda_s^+}(0,h+s)}t,
\end{equation}
and the lower correction at each integer $h\ge1$ satisfies
\begin{equation}\label{eq:boundary-lower-limit}
 \delta_h^-(t,x)
 =x\lim_{s\downarrow0}
       \frac{K_{\boldsymbol{\varphi},\lambda_s^-}(0,h-s)}t.
\end{equation}
In \eqref{eq:boundary-upper-limit} and \eqref{eq:boundary-lower-limit},
the limit as $s\downarrow0$ is taken separately for each coefficient
in the expansion of the kernel quotient in powers of $t$.
\end{lemma}

\begin{proof}
Fix $x\in[0,1]$ and $0<s<1$. On $0\le u<1$, define
\[
 \lambda_s^+(u)\coloneqq
 \begin{cases}
  x,&u=0,\\
  u/s,&0<u<s,\\
  1,&u=s,\\
  (u-s)/(1-s),&s<u<1,
 \end{cases}
 \qquad
 \lambda_s^-(u)\coloneqq
 \begin{cases}
  0,&u=0,\\
  u/(1-s),&0<u<1-s,\\
  x,&u=1-s,\\
  (u-1+s)/s,&1-s<u<1.
 \end{cases}
\]
Extend $\lambda_s^+$ and $\lambda_s^-$ with period $1$.
An affine bijection from an interval of length $L$ onto $(0,1)$ pushes
Lebesgue measure forward to $L\eta$. For each function the two open
subintervals have lengths $s$ and $1-s$, so the occupation measure on
a full period is $\eta$.

For integers $0\le j\le h$, there are $h-j$ full periods from $j$ to $h$.
On the remaining interval $(h,h+s)$, the occupation measure of
$\lambda_s^+$ is $s\eta$ by translation from $(0,s)$. Hence on
$(j,h+s)$ the occupation measure is $(h-j+s)\eta$. With the prescribed values
$\lambda_s^+(0)=x$ and $\lambda_s^+(h+s)=\lambda_s^+(s)=1$,
the kernel formula \eqref{eq:2-14} gives
\[
 \frac{K_{\boldsymbol{\varphi},\lambda_s^+}(0,h+s)}t
 =\sum_{j\le h}(\mathsf T_{h-j+s}\varphi_j)[x,1].
\]
For integers $0\le j<h$, there are $h-j-1$ full periods from $j$ to $h-1$.
On the remaining interval $(h-1,h-s)$, the occupation measure of
$\lambda_s^-$ is $(1-s)\eta$ by translation from $(0,1-s)$. Hence on
$(j,h-s)$ the occupation measure is $(h-j-1)\eta+(1-s)\eta=(h-j-s)\eta$.
With $\lambda_s^-(0)=0$ and $\lambda_s^-(h-s)=\lambda_s^-(1-s)=x$,
\eqref{eq:2-14} gives
\[
 \frac{K_{\boldsymbol{\varphi},\lambda_s^-}(0,h-s)}t
 =\sum_{j<h}(\mathsf T_{h-j-s}\varphi_j)[0,x].
\]
For fixed $x$ and each integer $k\ge0$, \eqref{eq:uniform-transport-monomial}
expresses the coefficient of $t^k$ in each kernel quotient as a polynomial
in $s$ on $0<s<1$. To take the limit of the coefficient as $s\downarrow0$,
evaluate that polynomial at $s=0$. The transport lengths $h-j+s$ and
$h-j-s$ then both reduce to $h-j$. By \eqref{eq:2-10}, multiplying the upper sum by $1-x$ gives
$\delta_h^+(t,x)$, and multiplying the lower sum by $x$ gives
$\delta_h^-(t,x)$.
\end{proof}

For a finite normal form $\boldsymbol{\varphi}$ representing
$\beta(t,z)$, an auxiliary function $\lambda$, and $u\in\R$,
define the measure in the displacement variable $\Delta=v-u$ by
\[
 \nu_{\boldsymbol{\varphi},\lambda,u}
 \coloneqq\sum_{j\ge0}\varphi_j(t,\lambda(u))\operatorname{Dirac}(j)
 +K_{\boldsymbol{\varphi},\lambda}(u,u+\Delta)\,\mathrm{d}\Delta.
\]
The point mass $\operatorname{Dirac}(j)$ from Definition~\ref{def:canonical-kernel}
records evaluation of the test function at $u+j$.
For each bounded Borel set $B\subseteq[0,+\infty)$, integration over $B$
coefficientwise in $t$ gives a polynomial $\nu_{\boldsymbol{\varphi},\lambda,u}(B)$
whose coefficients may have either sign.

Fix integers $0\le h_-<h_+$. For every positive integer $t$, the two
orders of closed truncation give
\[
 \bigl(\beta^{\ge h_-}\bigr)^{\le h_+}(t,z)
 =\bigl(\beta^{\le h_+}\bigr)^{\ge h_-}(t,z)
 =\sum_{a=th_-}^{th_+}([z^a]\beta(t,z))z^a.
\]
By Lemma~\ref{lem:closed-clipping}, truncating $\boldsymbol{\varphi}$ in
either order gives a normal form of $\left.\beta(t,z)\right|_{[th_-,\,th_+]}$.
Uniqueness in Lemma~\ref{lem:normal-form-uniqueness} therefore gives
\[
 \widetilde{\boldsymbol{\varphi}}
 \coloneqq\bigl(\boldsymbol{\varphi}^{\ge h_-}\bigr)^{\le h_+}
 =\bigl(\boldsymbol{\varphi}^{\le h_+}\bigr)^{\ge h_-}.
\]

To compute either endpoint component, truncate $\boldsymbol{\varphi}$ first
at that endpoint. Since $h_-<h_+$, the upper truncation preserves the
component at $h_-$, and the lower truncation preserves the component at $h_+$.
Thus \eqref{eq:2-10} applied to the original $\boldsymbol{\varphi}$ gives
\[
 \widetilde{\varphi}_{h_-}=\varphi_{h_-}+\delta_{h_-}^-,
 \qquad
 \widetilde{\varphi}_{h_+}=\varphi_{h_+}+\delta_{h_+}^+.
\]
For $h_-<j<h_+$, $\widetilde{\varphi}_j=\varphi_j$,
and for $j\notin[h_-,h_+]$, $\widetilde{\varphi}_j=0$.

Applying Lemma~\ref{lem:closed-clipping} first at the lower bound $h_-$
and then at the upper bound $h_+$ gives, for every positive noninteger $\Delta$,
\begin{equation}\label{eq:double-cut-density}
 K_{\widetilde{\boldsymbol{\varphi}},\lambda}(u,u+\Delta)
 =\mathbf1_{\{h_-<\Delta<h_+\}}
       K_{\boldsymbol{\varphi},\lambda}(u,u+\Delta).
\end{equation}
Combining the point mass of weight $\widetilde{\varphi}_j(t,\lambda(u))$
at each integer $j\ge0$ with the density in \eqref{eq:double-cut-density} gives
\begin{equation}\label{eq:complete-boundary-measure}
\begin{aligned}
 \nu_{\widetilde{\boldsymbol{\varphi}},\lambda,u}
 &=\mathbf1_{[h_-,h_+]}\nu_{\boldsymbol{\varphi},\lambda,u}\\
 &\quad+\delta_{h_-}^-(t,\lambda(u))\operatorname{Dirac}(h_-)
       +\delta_{h_+}^+(t,\lambda(u))\operatorname{Dirac}(h_+).
\end{aligned}
\end{equation}
For each bounded Borel set $B\subseteq[0,+\infty)$, the restriction
$\mathbf1_{[h_-,h_+]}\nu_{\boldsymbol{\varphi},\lambda,u}$ assigns to $B$
the value $\nu_{\boldsymbol{\varphi},\lambda,u}(B\cap[h_-,h_+])$.

\subsection{Positive states and interval operators}\label{subsec:positive-states}

To prove Ehrhart positivity, we will show that for every nonempty interval
$I\subseteq N$, the evaluation weights and the integral kernel of the operator
representing $\beta_I(t,z)$ have nonnegative coefficients in $t$.
We first verify nonnegativity for a single coordinate. For an interval
partition $\pi$ of $I$, the evaluation weight at an integer jump $h$ for
$\prod_{J\in\pi}\beta_J(t,z)$ is the sum of the products of block evaluation
weights over all choices of integer block jumps totaling $h$. By
Lemma~\ref{lem:mixed-kernel-composition}, the product kernel is obtained by
summing and integrating products of block weights and kernels. Thus
nonnegative coefficients in $t$ are preserved under block products.

\begin{definition}\label{def:complete-positive-state}
Fix nonnegative integers $h_-\le h_+$.
A finite normal form $\boldsymbol{\varphi}=(\varphi_j)_{j\ge0}$ of a lattice-point
generating function $\beta(t,z)$ is a \emph{complete positive state} of
support $[h_-,h_+]$ if the following conditions hold for every
auxiliary function $\lambda$ and all $j\in\Z_{\ge0}$, $u,v\in\R$:
\begin{enumerate}[label=(\roman*),ref=(\roman*),leftmargin=*,topsep=4pt,itemsep=3pt,parsep=0pt]
\item\label{cond:state-positivity}
The evaluation weight $\varphi_j(t,\lambda(u))$ and the integral
kernel $K_{\boldsymbol{\varphi},\lambda}(u,v)$ have nonnegative
coefficients in $t$.
\item\label{cond:state-support}
The weight $\varphi_j(t,\lambda(u))$ vanishes when
$j\notin[h_-,h_+]$, and $K_{\boldsymbol{\varphi},\lambda}(u,v)$
vanishes when $v-u\notin[h_-,h_+]$.
\item\label{cond:state-endpoints}
The components satisfy
\begin{equation}\label{eq:2-12}
 \varphi_j(t,1)=0\quad\text{if }j\ne h_-,\qquad
 \varphi_j(t,0)=0\quad\text{if }j\ne h_+.
\end{equation}
\end{enumerate}
\end{definition}

\begin{example}\label{ex:running-singleton}
The projection of $\mathcal P_\star$ onto its first coordinate is $[0,2]$.
For any integer $c\ge1$, the lattice-point generating function of $[0,c]$
has finite normal form $\boldsymbol\alpha_c$. By \eqref{eq:2-13},
for $g\in\mathscr F[t]$,
\[
 (\mathsf R_{\lambda,\boldsymbol\alpha_c}g)(u)
 =\lambda(u)g(u)+(1-\lambda(u))g(u+c)
   +t\int_0^c g(u+s)\,\mathrm{d}s.
\]
The only evaluation weights are $\lambda(u)$ at displacement $0$ and
$1-\lambda(u)$ at displacement $c$, both nonnegative.
The integral kernel $K_{\boldsymbol\alpha_c,\lambda}(u,u+\Delta)$
equals $t$ for noninteger $\Delta\in(0,c)$ and zero otherwise.
Thus the weights and kernel have nonnegative coefficients in $t$
and vanish outside $[0,c]$.
When $\lambda(u)=0$, the sole evaluation term is $g(u+c)$.
When $\lambda(u)=1$, it is $g(u)$. The integral term is independent
of $\lambda(u)$ and remains in both cases.
Hence $\boldsymbol\alpha_c$ is a complete positive state on $[0,c]$.
\end{example}

\begin{lemma}\label{lem:positive-products}
If $\boldsymbol{\varphi}$ and $\boldsymbol{\psi}$ are complete positive states on
$[h_-,h_+]$ and $[\widetilde h_-,\widetilde h_+]$, respectively,
then $\boldsymbol{\varphi}\boldsymbol{\psi}$ is a complete positive state on
$[h_-+\widetilde h_-,h_++\widetilde h_+]$.
\end{lemma}

\begin{proof}
The convolution $\boldsymbol{\varphi}\boldsymbol{\psi}$ represents the product
of the two generating functions. For each integer $h\ge0$, the convolution formula gives
\[
 (\boldsymbol{\varphi}\boldsymbol{\psi})_h(t,x)
 =\sum_{a+b=h}\varphi_a(t,x)\psi_b(t,x).
\]
Fix an auxiliary function $\lambda$ and $u\in\R$. At $x=\lambda(u)$,
each factor in the sum has nonnegative coefficients in $t$, and a nonzero
product requires $a\in[h_-,h_+]$ and $b\in[\widetilde h_-,\widetilde h_+]$.
Hence the evaluation weight has nonnegative coefficients and vanishes unless
$h$ lies in $[h_-+\widetilde h_-,h_++\widetilde h_+]$.
By \eqref{eq:2-12}, a nonzero product at $x=1$ requires
$(a,b)=(h_-,\widetilde h_-)$, and at $x=0$ it requires
$(a,b)=(h_+,\widetilde h_+)$. Thus the only possible indices $h$ at $x=1$
and $x=0$ are $h_-+\widetilde h_-$ and $h_++\widetilde h_+$, respectively.

For positive noninteger $\Delta=v-u$, each summand in the two sums
in \eqref{eq:2-16} is an evaluation weight times a kernel value, and
the integrand is a product of two kernel values. All factors have nonnegative coefficients in $t$
by condition~\ref{cond:state-positivity}. Summation and coefficientwise
integration over $w\in[u,v]$ therefore give nonnegative coefficients for
$K_{\boldsymbol{\varphi}\boldsymbol{\psi},\lambda}(u,v)$.
The displacements for $\boldsymbol{\varphi}$ and $\boldsymbol{\psi}$ in
the first sum, the second sum, and the integrand are, respectively,
\[
 (a,\Delta-a),\qquad (\Delta-b,b),\qquad (w-u,v-w).
\]
In each pair the displacements add to $\Delta$. By
condition~\ref{cond:state-support}, a nonzero product requires the first
displacement to lie in $[h_-,h_+]$ and the second in
$[\widetilde h_-,\widetilde h_+]$. Thus the product kernel vanishes outside
$[h_-+\widetilde h_-,h_++\widetilde h_+]$. At integer or nonpositive
displacements the kernel is zero by Definition~\ref{def:canonical-kernel}.
\end{proof}

We construct a finite normal form for $\beta_I(t,z)$ by induction on the
length of the nonempty interval $I\subseteq N$. For each $i\in N$, set
$\boldsymbol{\varphi}_{\{i\}}\coloneqq\boldsymbol{\alpha}_{c_i}$,
which represents $\beta_{\{i\}}(t,z)$ by \eqref{eq:singleton-generating-function}
and Example~\ref{ex:running-singleton} with $c=c_i$.
For $|I|\ge2$, every block $J$ of a partition of $I$ into at least two
intervals is shorter than $I$, so $\boldsymbol{\varphi}_J$ is already defined.
The integer bounds satisfy $0\le\underline b(I)<\overline b(I)$ by
normalization and \eqref{eq:1-5}. Using convolution for the block products, set
\begin{equation}\label{eq:4-1}
\begin{aligned}
 \boldsymbol{\psi}_I&\coloneqq\sum_{\substack{\pi\in\operatorname{Int}(I)\\|\pi|\ge2}}
       (-1)^{|\pi|}\prod_{J\in\pi}\boldsymbol{\varphi}_J,\\
 \boldsymbol{\varphi}_I&\coloneqq\bigl(\boldsymbol{\psi}_I^{\ge \underline b(I)}\bigr)^{\le \overline b(I)}.
\end{aligned}
\end{equation}
Write $\boldsymbol{\varphi}_I=(\varphi_{I,j}(t,x))_{j\ge0}$ and, for
$|I|\ge2$, $\boldsymbol{\psi}_I=(\psi_{I,j}(t,x))_{j\ge0}$.

Suppose that $\boldsymbol{\varphi}_J$ represents $\beta_J(t,z)$ for every
proper subinterval $J$ of $I$, where $|I|\ge2$. By Section~\ref{subsec:normal-forms},
each block convolution in \eqref{eq:4-1} represents
$\prod_{J\in\pi}\beta_J(t,z)$, and summing with signs $(-1)^{|\pi|}$ gives
$\gamma_I(t,z)$.
For each positive integer $t$, put $x=(1-z)^{-1}$. By
Section~\ref{subsec:closed-truncation}, truncating $\boldsymbol{\psi}_I$
as in \eqref{eq:4-1} retains precisely the terms of $\gamma_I(t,z)$ with
exponents from $t\underline b(I)$ through $t\overline b(I)$. Hence
\[
\begin{aligned}
 \sum_{j\ge0}z^{jt}\psi_{I,j}(t,x)&=\gamma_I(t,z),\\
 \sum_{j\ge0}z^{jt}\varphi_{I,j}(t,x)
 &=\left.\gamma_I(t,z)\right|_{[t\underline b(I),\,t\overline b(I)]}
 =\beta_I(t,z),
\end{aligned}
\]
where the last equality is the counting recursion \eqref{eq:counting-recursion}.

Fix an auxiliary function $\lambda$. For $|I|\ge2$, list the blocks of each
$\pi\in\operatorname{Int}(I)$ as $H_1,\ldots,H_s$ in the linear order on $I$.
Formula \eqref{eq:normal-form-operator} makes the map
$\boldsymbol{\varphi}\mapsto\mathsf R_{\lambda,\boldsymbol{\varphi}}$
linear, and Lemma~\ref{lem:mixed-kernel-composition} identifies the operator
of each block convolution with the composition of the block operators.
To compute the evaluation weights and integral kernel of
$\mathsf R_{\lambda,\boldsymbol\psi_I}$ block by block in
Lemma~\ref{lem:exact-inventory}\ref{part:inventory-before}, we apply this map
to the signed sum in \eqref{eq:4-1} and obtain
\begin{equation}\label{eq:interval-operator-expansion}
 \mathsf R_{\lambda,\boldsymbol\psi_I}
 =\sum_{\substack{\pi\in\operatorname{Int}(I)\\|\pi|\ge2}}(-1)^s
   \mathsf R_{\lambda,\boldsymbol\varphi_{H_1}}\cdots
   \mathsf R_{\lambda,\boldsymbol\varphi_{H_s}}.
\end{equation}

\section{Cut cancellation and nonnegative truncation corrections}\label{sec:recursion}

Lemma~\ref{lem:positive-products} shows that convolution preserves complete
positive states. The interval recursion \eqref{eq:4-1} also takes the
signed sum $\boldsymbol{\psi}_I$ of the block products and truncates it
at $\underline b(I)$ and $\overline b(I)$. Proposition~\ref{prop:state-closure},
proved at the beginning of Section~\ref{sec:main-proof}, gives the positivity,
support, and endpoint conditions for the resulting normal form
$\boldsymbol{\varphi}_I$.

\begin{proposition}\label{prop:state-closure}
For every nonempty interval $I\subseteq N$,
$\boldsymbol{\varphi}_I$ is a complete positive state on $[\underline b(I),\overline b(I)]$.
\end{proposition}

Section~\ref{subsec:sources} computes the signed multiplicity $\sigma(y)$
of each source in $\boldsymbol{\psi}_I$. For total jumps
$y(I)\in[\underline b(I),\overline b(I)]$, this multiplicity is zero or one.
For $0<s<1$, Section~\ref{sec:boundary} bounds the measure of the region
where $\sigma(y)<0$ on each of the two fibers $y(I)=\underline b(I)-s$
and $y(I)=\overline b(I)+s$. The fiber estimate and the one-sided limits
as $s\downarrow0$ are used to prove that both endpoint corrections have
nonnegative coefficients in $t$.
Section~\ref{sec:main-proof} then uses Proposition~\ref{prop:state-closure}
to recover $i(\mathcal P,t)$ from $[z^{rt}](\beta_N\beta_{\{e\}})$.

\subsection{Weighted cancellation over interval partitions}\label{subsec:sources}

Fix an interval $I\subseteq N$ with $|I|\ge2$, and put
$h_-\coloneqq\underline b(I)$ and $h_+\coloneqq\overline b(I)$.
For each interval $J\subseteq I$ with $|J|\ge2$, let $\delta_J^-$ and
$\delta_J^+$ denote the lower and upper corrections computed from
$\boldsymbol{\psi}_J$ by \eqref{eq:2-10} at $\underline b(J)$ and
$\overline b(J)$, respectively. Applying the endpoint component identities
of Section~\ref{subsec:closed-truncation} to $\boldsymbol{\psi}_I$
with cutoffs $h_-,h_+$ gives
\begin{equation}\label{eq:interval-endpoint-components}
 \varphi_{I,h_-}=\psi_{I,h_-}+\delta_I^-,
 \qquad
 \varphi_{I,h_+}=\psi_{I,h_+}+\delta_I^+.
\end{equation}

For nonempty intervals $A,B\subseteq I$, write
$A\prec B$ if every element of $A$ precedes every element of $B$
in the fixed linear order on $I$. Fix an auxiliary function $\lambda$
and an initial position $u\in\R$.

\begin{definition}\label{def:sources}
A \emph{source} on $I$ consists of an interval partition
$I=K_1\sqcup\cdots\sqcup K_q$, with the blocks listed in the
linear order on $I$, together with a branch choice on each block.
Each block is a singleton $\{i\}$ or a
\emph{correction block} $J$ with $|J|\ge2$.
For each block $J$, the chosen branch increases the argument of the test
function by its \emph{jump} $y_J$. The \emph{entry position} of $J$ is
\[
 v\coloneqq u+\sum_{J'\prec J}y_{J'},
\]
where $J'$ ranges over the source blocks preceding $J$.
The argument $v+y_J$ is passed to the next block, or to the test
function if $J=K_q$.
When block $J$ acts on a test function $g$, the \emph{weight} of an atomic
branch is the factor multiplying $g(v+y_J)$ in the evaluation term.
For a continuous branch, the weight is the factor multiplying
$g(v+y_J)\,\mathrm{d}y_J$ in the integral.
For a singleton block $\{i\}$ entered at $v$, the pairs of jump and
weight for the two atomic branches in \eqref{eq:2-13} are
\[
 \bigl(0,\,\lambda(v)\bigr),
 \qquad
 \bigl(c_i,\,1-\lambda(v)\bigr).
\]
On the continuous branch, we integrate the jump $y_{\{i\}}$ over
$(0,c_i)$ with weight $t$ and measure $\mathrm{d}y_{\{i\}}$.
A correction block $J$ entered at $v$ chooses one of the two endpoint
correction terms in \eqref{eq:complete-boundary-measure} when $\boldsymbol{\psi}_J$
is truncated at $\underline b(J)$ and $\overline b(J)$.
For the labels $\delta_J^-$ and $\delta_J^+$, respectively, the pairs
of jump and weight are
\[
 \bigl(\underline b(J),\,\delta_J^-(t,\lambda(v))\bigr),
 \qquad
 \bigl(\overline b(J),\,\delta_J^+(t,\lambda(v))\bigr).
\]

Given the block jumps, define the auxiliary vector $y=(y_j)_{j\in I}$ by
$y_j\coloneqq y_{K_i}/|K_i|$ for $j\in K_i$. When $q\ge2$, the partition
$K_1,\ldots,K_q$ and $y$ determine the signed multiplicity $\sigma(y)$
in \eqref{eq:4-3}.
For a nonempty union $H$ of consecutive whole blocks, the \emph{guard}
is the closed constraint
\[
 \underline b(H)\le y(H)=\sum_{K_i\subseteq H}y_{K_i}\le\overline b(H).
\]
For a correction block $J$, the branch specifies only the total jump $y_J$
and imposes no guard on a proper subinterval of $J$.
\end{definition}

\begin{example}\label{ex:running-source}
Take $I=N$ with the four singleton blocks from Example~\ref{ex:running-cuts}.
Choose the continuous branches on coordinates $1$ and $3$, the upper
endpoint branch on coordinate $2$, and the lower endpoint branch on
coordinate $4$.
\[
\begin{array}{c|c|c|c}
 \text{block}&\text{jump}&\text{entry position}&\text{weight}\\ \hline
 \{1\}&a&u&t\\
 \{2\}&2&u+a&1-\lambda(u+a)\\
 \{3\}&b&u+a+2&t\\
 \{4\}&0&u+a+2+b&\lambda(u+a+2+b)
\end{array}
\]
The continuous jumps $a$ and $b$ range independently over $(0,2)$ and
$(0,1)$ with product measure $\mathrm{d}a\,\mathrm{d}b$. Jumps of $0$ or $2$
on coordinate $1$ and of $0$ or $1$ on coordinate $3$ belong to the atomic branches.
Imposing only the guard $2\le y(I)=2+a+b\le3$ gives the vectors
$y=(a,2,b,0)$ of Example~\ref{ex:running-overlap} with $a,b>0$ and
$a+b\le1$, with equality when $y(I)=3$.
\end{example}

\begin{definition}
\label{def:candidate-source-inventory}
For $I=\{i\}\subseteq N$, the \emph{candidate sources} for
$\boldsymbol{\varphi}_I$ are the three branches of
$\mathsf R_{\lambda,\boldsymbol\alpha_{c_i}}$ in \eqref{eq:2-13}.
For an interval $I\subseteq N$ with $|I|\ge2$, the candidate sources
for $\boldsymbol{\varphi}_I$ are as follows.
\begin{enumerate}[label=(\roman*),leftmargin=*]
\item For each partition of $I$ into consecutive blocks $K_1,\ldots,K_q$
with $q\ge2$, include one source for each choice of a branch on each block.
Impose $\underline b(H)\le y(H)\le\overline b(H)$ for every
$H=K_a\cup\cdots\cup K_d$ with $1\le a\le d\le q$, including $H=I$.
\item Include the two single-block sources labelled $\delta_I^-$ and
$\delta_I^+$ once each, at jumps $\underline b(I)$ and $\overline b(I)$,
respectively. Neither source imposes a constraint on a proper subinterval of $I$.
\end{enumerate}
\end{definition}

Fix a source on $I$ with blocks $K_1,\ldots,K_q$. Let $D$ consist of the
indices $i$ for which $\{i\}$ uses the continuous branch. The jumps
$y_i\in(0,c_i)$ for $i\in D$ are the integration variables. All other
blocks have fixed integer jumps with sum $h$, so $y(I)=\Delta=h+\sum_{i\in D}y_i$.
For the fixed $\lambda$ and $u$, let $\upsilon(t,y)$ be the product
of the branch weights at the entry positions determined by $y$.
For $q\ge2$, let $\mathcal B$ be the family of block unions $H$ in
Definition~\ref{def:candidate-source-inventory}(i). For $q=1$, set
$\mathcal B=\varnothing$, since the branch range already enforces the bounds on $y(I)$.

If $D=\varnothing$, all jumps are fixed and $\Delta=h$.
For an input $g\in\mathscr F[t]$, the source contributes
the evaluation term
\[
 \upsilon(t,y)
 \prod_{H\in\mathcal B}\mathbf1_{\{\underline b(H)\le y(H)\le \overline b(H)\}}
 \,g(u+h).
\]
If $D\ne\varnothing$, choose $j\in D$. The change of variables
\begin{equation}\label{eq:source-fiber-coordinates}
 (y_i)_{i\in D}\longmapsto
 \bigl(\Delta,(y_i)_{i\in D\setminus\{j\}}\bigr)
\end{equation}
has Jacobian of absolute value one. On the fiber $y(I)=\Delta$, write
$y_j=\Delta-h-\sum_{i\in D\setminus\{j\}}y_i$ and let $\mu_\Delta$ be
product Lebesgue measure in the remaining coordinates, restricted by
$0<y_i<c_i$ for every $i\in D$. For $k\in D\setminus\{j\}$, fixing
$\Delta$ and all $y_i$ with $i\in D\setminus\{j,k\}$ fixes $y_j+y_k$.
Exchanging $y_k$ for $y_j$ as a free coordinate therefore has slope $-1$
and preserves $\mu_\Delta$.
If $D=\{j\}$, $\mu_\Delta$ assigns mass one to $y_j=\Delta-h$ when
$0<\Delta-h<c_j$ and is zero otherwise. By \eqref{eq:source-fiber-coordinates},
the density of the source with respect to $\mathrm{d}\Delta$ is
\[
 \int_{y(I)=\Delta}\upsilon(t,y)
 \prod_{H\in\mathcal B}\mathbf1_{\{\underline b(H)\le y(H)\le\overline b(H)\}}
 \,\mathrm{d}\mu_\Delta(y).
\]
Thus a source with a continuous branch gives no point mass in $\Delta$,
even when $|D|=1$.

Fix a positive noninteger $\Delta$ and let $H$ be a union of consecutive
source blocks. On the fiber $y(I)=\Delta$, we claim that
\[
 H\cap D\ne\varnothing \quad\Longrightarrow\quad
 \mu_\Delta\{y:y(H)\in\Z\}=0.
\]
If $D\subseteq H$, then $y(H)=\Delta-y(I\setminus H)$ is nonintegral
because all jumps outside $H$ are integers. If
$\varnothing\ne H\cap D\subsetneq D$, choose
$i\in H\cap D$ and $j\in D\setminus H$ and eliminate $y_j$ in
\eqref{eq:source-fiber-coordinates}. With $y_k$ fixed for
$k\in D\setminus\{i,j\}$, the sum $y(H)$ has slope one in $y_i$.
Only countably many values of $y_i$ can then give $y(H)\in\Z$, so
integrating $\mathbf1_{\{y(H)\in\Z\}}$ first in $y_i$ gives zero.
If $H\cap D=\varnothing$, the sum $y(H)$ is fixed. When it equals either
$\underline b(H)$ or $\overline b(H)$, the indicator of the closed
constraint on $H$ is one throughout the fiber.

For a partition $\pi$ in the sum defining $\boldsymbol\psi_I$ in
\eqref{eq:4-1}, expand $\prod_{H\in\pi}\boldsymbol\varphi_H$ by choosing
one candidate source for each $H\in\pi$ and concatenating the sources
in block order. If $K_1,\ldots,K_q$ are the blocks of the concatenated
source, each factor $\boldsymbol\varphi_H$ tests only the bounds on
unions of consecutive $K_i$ contained in $H$.
As in Section~\ref{sec:violations}, $\mathcal H$ consists of all nonempty
proper unions of consecutive $K_i$, and $\mathcal C_y$ consists of the
nonempty subsets of gaps that cut every violating member of $\mathcal H$.

\begin{samepage}
\begin{lemma}\label{lem:weighted-cut-cancellation}
Fix a source $\mathfrak s$ with at least two blocks. For each choice of
jumps $y$ in its branch ranges, the terms in the candidate expansion
of $\boldsymbol\psi_I$ whose block sources concatenate to $\mathfrak s$
have total signed weight $\sigma(y)\upsilon(t,y)$.
\end{lemma}
\end{samepage}

\begin{proof}
Write $K_1,\ldots,K_q$ for the blocks of $\mathfrak s$. A partition
in \eqref{eq:4-1} that splits a $K_i$ cannot produce $\mathfrak s$ by
concatenation. The remaining partitions are $\pi=\pi(C)$ for nonempty
subsets $C$ of the gaps between consecutive $K_i$. If $H\in\pi$ is a
single $K_i$, keep its branch, including the label on a correction
block. The branch range ensures that the guard on $K_i$ holds.
If $H$ contains several $K_i$, take them in order with their prescribed
branches. Definition~\ref{def:candidate-source-inventory}(i) then imposes
all guards on nonempty unions of consecutive $K_i$ contained in $H$.
These restrictions of $\mathfrak s$ give the unique choice on $\pi$
that concatenates to $\mathfrak s$, even when a weight or guard
indicator vanishes.

For $H\in\pi$ and a source block $J\subseteq H$, first entering $H$
and then advancing to $J$ gives the entry position
\[
 u+\sum_{J'\prec H}y_{J'}
   +\sum_{\substack{J'\subseteq H\\J'\prec J}}y_{J'}
 =u+\sum_{J'\prec J}y_{J'},
\]
where $J'$ ranges over the blocks of $\mathfrak s$. Each branch thus
uses the same value of $\lambda$ for every $\pi$, giving the common
weight product $\upsilon(t,y)$. For $D=\varnothing$ the evaluation is
at $u+h$. For $D\ne\varnothing$ the continuous branches of $\mathfrak s$
fix the integration variables and their ranges, so
\eqref{eq:source-fiber-coordinates} gives the same $\mu_\Delta$ for every
$\pi$, with the closed guard indicators in the integrand.

A union $G\in\mathcal H$ lies in one block of $\pi(C)$ if and only if
$C$ contains no gap inside $G$. Thus the product of the guard indicators
within the blocks of $\pi(C)$ equals $\mathbf1_{\{C\in\mathcal C_y\}}$.
Using $|\pi(C)|=|C|+1$ and \eqref{eq:4-3}, we obtain
\begin{equation}\label{eq:weighted-cut-cancellation}
 \sum_{\pi}(-1)^{|\pi|}\upsilon(t,y)
   \prod_{H\in\pi}\prod_{\substack{G\in\mathcal H\\G\subseteq H}}
     \mathbf1_{\{\underline b(G)\le y(G)\le \overline b(G)\}}
 =\upsilon(t,y)\sum_{C\in\mathcal C_y}(-1)^{|C|+1}
 =\sigma(y)\upsilon(t,y)
\end{equation}
as an identity in $\R[t]$.
\end{proof}

By Lemma~\ref{lem:cut-sum}, applied to the fixed blocks $K_i$
and jumps $y$, if $\underline b(I)\le y(I)\le \overline b(I)$, the combined signed
weight in \eqref{eq:weighted-cut-cancellation} is $\upsilon(t,y)$
when every union of consecutive $K_i$ satisfies its guard,
and zero otherwise.

\begin{example}\label{ex:running-transport}
Take the endpoint jumps $y_2=2$, $y_4=0$ and the continuous jumps
$y_1=a$, $y_3=b$ of Example~\ref{ex:running-source}.
Impose the total guard $a+b\le1$. The two partitions allowed by the
block guards in Example~\ref{ex:running-cuts} differ only in whether
coordinates $3$ and $4$ share a block.
Figure~\ref{fig:running-source-cancellation} shows where the jump $b$
enters the argument at coordinate $4$ in each partition.
Grouping coordinates $3$ and $4$ preserves the range of $b$, since
$(b,0)\in\mathcal P_{\{3,4\}}=[0,1]^2$ for every $0<b<1$.
For $g\in\mathscr F[t]$, the unsigned term for either partition is
\[
 t^2\!\int_{\substack{a>0,\ b>0\\a+b\le1}}
 (1-\lambda(u+a))\lambda(u+a+2+b)
 g(u+a+2+b)\,\mathrm{d}a\,\mathrm{d}b.
\]
The opposite signs give pointwise cancellation, including at $a+b=1$,
where $y(I)=3$.
\end{example}

\begin{figure}[htbp]
\centering
\begin{tikzpicture}[font=\small]
 \matrix[column sep=8mm,row sep=2.5mm,nodes={anchor=base,inner sep=0pt}] {
   \node {Partition $\pi$}; & \node {Sign}; &
   \node {Entry at coordinate $4$}; \\
   \node {$\{1\}\mid\{2\}\mid\{3,4\}$}; & \node {$-1$}; &
   \node {$\displaystyle
     \underbrace{u+a+2}_{\text{entry into }\{3,4\}}
     +\underbrace{b}_{\text{jump at }3}$}; \\
   \node {$\{1\}\mid\{2\}\mid\{3\}\mid\{4\}$}; & \node {$+1$}; &
   \node {$\displaystyle
     \underbrace{u+a+2+b}_{\text{entry into }\{4\}}$}; \\
 };
\end{tikzpicture}
\caption{In both partitions of Example~\ref{ex:running-transport}, the
zero jump at coordinate $4$ has weight $\lambda(u+a+2+b)$.}
\label{fig:running-source-cancellation}
\end{figure}
\FloatBarrier

\Needspace{12\baselineskip}
\begin{lemma}\label{lem:exact-inventory}
Let $I\subseteq N$ be a nonempty interval, $\lambda$ an auxiliary
function, and $u\in\R$. Fix an integer $j\ge0$ and a noninteger $\Delta>0$.
\begin{enumerate}[label=(\roman*),ref=(\roman*),leftmargin=*]
\item\label{part:inventory-before} For $|I|\ge2$, summing over sources with at least two blocks
before imposing $\underline b(I)\le y(I)\le\overline b(I)$ gives
\[
\begin{aligned}
 \psi_{I,j}(t,\lambda(u))
 &=\sum_{\substack{D=\varnothing\\y(I)=j}}\sigma(y)\upsilon(t,y),\\
 K_{\boldsymbol{\psi}_I,\lambda}(u,u+\Delta)
 &=\sum_{D\ne\varnothing}\int_{y(I)=\Delta}
     \sigma(y)\upsilon(t,y)\,\mathrm{d}\mu_{\Delta}(y).
\end{aligned}
\]
\item\label{part:inventory-after} For the candidate sources of Definition~\ref{def:candidate-source-inventory},
the weights for $D=\varnothing$ and $y(I)=j$ sum to $\varphi_{I,j}(t,\lambda(u))$,
and the fiber integrals at $y(I)=\Delta$ for $D\ne\varnothing$ sum to
$K_{\boldsymbol{\varphi}_I,\lambda}(u,u+\Delta)$.
\end{enumerate}
\end{lemma}

\begin{proof}
For $I=\{i\}$, \ref{part:inventory-after} follows from \eqref{eq:2-13}.
Assume \ref{part:inventory-after} on every proper subinterval of $I$ at
every initial position. Fix a partition $\pi$ in \eqref{eq:interval-operator-expansion}
and one candidate source on each block $H\in\pi$.
In the ordered block product, the coefficient of $g(u+j)$ is the sum of
$\upsilon(t,y)$ times the within-block guard indicators over the choices
with $D=\varnothing$ and $y(I)=j$, by Lemma~\ref{lem:mixed-kernel-composition}.

For $D\ne\varnothing$, work on the fiber $y(I)=\Delta$. Split the product
before the last block, with coordinate sets $A$ on the left and $B$ on the right.
If both sets meet $D$, choose $i\in A\cap D$ and $\ell\in B\cap D$.
Eliminate $y_\ell$ by \eqref{eq:source-fiber-coordinates} and fix $y_r$
for $r\in D\setminus\{i,\ell\}$. Since $\ell\notin A$, the sum
$y(A\setminus\{i\})$ is fixed in the change of variable
\[
 a=y(A)=y_i+y(A\setminus\{i\}),\qquad
 \mathrm{d}a=\mathrm{d}y_i,\qquad y(B)=\Delta-a .
\]
At fixed $a$, the remaining coordinates carry the product of the measures
in \eqref{eq:source-fiber-coordinates} on $y(A)=a$ and $y(B)=\Delta-a$.
The right block starts at $w=u+a$ in the integral in \eqref{eq:2-16}.
If $A\cap D=\varnothing$ and $y(A)=h$, the first sum in \eqref{eq:2-16}
multiplies the right kernel at $(u+h,u+\Delta)$ by the product of atomic
weights on $A$ starting at $u$. If $B\cap D=\varnothing$ and $y(B)=h$,
the second sum multiplies the left kernel at $(u,u+\Delta-h)$ by the
product of atomic weights on $B$ starting at $u+\Delta-h$.

On the fixed fiber, $y(G)\notin\Z$ for $\mu_\Delta$-almost every $y$
simultaneously for all unions $G$ of consecutive blocks of $\pi$ with
$G\cap D\ne\varnothing$,
by the null-set calculation following \eqref{eq:source-fiber-coordinates}
and the finiteness of $\pi$. Repeat the last-block split by \eqref{eq:2-16}
while the left product contains a continuous branch and at least two blocks.
For a purely atomic prefix $G$, use the evaluation at $u+y(G)$.
For a single block $H\in\pi$, apply the induction hypothesis at
$v_H=u+\sum_{H'\in\pi,\,H'\prec H}y(H')$.
By Definition~\ref{def:sources}, each $[t^k]\upsilon(t,y)$ is a polynomial
in the values of the bounded Borel function $\lambda$ at the block entries.
The fiber has finite measure, so Fubini's theorem gives the term of
$K_{\prod_{H\in\pi}\boldsymbol\varphi_H,\lambda}(u,u+\Delta)$ for the
chosen block sources, coefficient by coefficient, as
\[
 \int_{y(I)=\Delta}\upsilon(t,y)
 \prod_{H\in\pi}\ \prod_{\substack{G\in\mathcal H\\G\subseteq H}}
 \mathbf1_{\{\underline b(G)\le y(G)\le\overline b(G)\}}\,
 \mathrm{d}\mu_\Delta(y).
\]
For a fixed concatenated source, \eqref{eq:weighted-cut-cancellation}
sums over the partitions producing that source with sign $(-1)^{|\pi|}$,
giving $\sigma(y)\upsilon(t,y)$ in both the atomic sum and the fiber
integral in \ref{part:inventory-before}.

On $h_-\le y(I)\le h_+$, Lemma~\ref{lem:cut-sum} makes $\sigma(y)$
the indicator of the guards in Definition~\ref{def:candidate-source-inventory}(i).
For each noninteger $\Delta>0$, \eqref{eq:double-cut-density} restricts
the kernel to $h_-\le\Delta\le h_+$ since $h_-,h_+\in\Z$.
For integer jumps, \eqref{eq:complete-boundary-measure} keeps the
evaluation weights in $[h_-,h_+]$ and adds each single-block weight
in Definition~\ref{def:candidate-source-inventory}(ii) once,
$\delta_I^-(t,\lambda(u))$ at $h_-$ and $\delta_I^+(t,\lambda(u))$ at $h_+$.
\end{proof}

\subsection{Nonnegativity of the truncation corrections}\label{sec:boundary}

Fix an interval $I\subseteq N$ with $|I|\ge2$ and write
$h_-=\underline b(I)$ and $h_+=\overline b(I)$.
Lemma~\ref{lem:cut-sum} gives $\sigma(y)\in\{0,1\}$ for $h_-\le y(I)\le h_+$.
To control the endpoint corrections $\delta_I^\pm$, we bound the measure
of the region where $\sigma(y)<0$ just outside each endpoint.

Let $\mathcal S_I$ be the set of sources on $I$ with at least two blocks
and at least one continuous branch. For each $\mathfrak s\in\mathcal S_I$,
$\mu_\Delta$ denotes the fiber measure of $\mathfrak s$ at $y(I)=\Delta$.

\begin{lemma}\label{lem:thin-negative-fibers}
There is a constant $c_I>0$, depending only on the interval bounds
$\underline b(H),\overline b(H)$ for $H\subseteq I$ and the lengths $c_i$
for $i\in I$, such that for $0<s<1$ and each $\Delta\in\{h_++s,h_--s\}$,
\[
 \sum_{\mathfrak s\in\mathcal S_I}
 \mu_\Delta\bigl(\{y:\sigma(y)<0\}\bigr)\le c_I s .
\]
\end{lemma}

\begin{proof}
Sources without continuous branches have integer total jump, so their
fibers are empty at both $h_++s$ and $h_--s$.
Fix $\mathfrak s\in\mathcal S_I$ with blocks $K_1,\ldots,K_q$,
and let $D$ index its continuous branches. For a union $A$ of consecutive blocks, put
\[
 \Delta_A=
 \begin{cases}
  y(A)-\overline b(A),&\Delta=h_++s,\\
  \underline b(A)-y(A),&\Delta=h_--s.
 \end{cases}
\]
In particular, $\Delta_I=s$ on the fiber $y(I)=\Delta$.

For any $y$ on this fiber with $\sigma(y)<0$, \eqref{eq:4-4} gives a
group $\mathcal J$ of minimal violating members of $\mathcal H$ with
$\bigcup_{H\in\mathcal J}H=I$. Each member is a proper subset of $I$, so
$\mathcal J$ has at least two members. In order of increasing left endpoint,
let $J$ be the last member and $V$ the union of all preceding members.
The proof of Lemma~\ref{lem:minimal-violations} shows that every member of
$\mathcal J$ and the union of every initial segment of $\mathcal J$ violate
the same bound as $I$. In particular, $\Delta_J>0$ and $\Delta_V>0$.
The member preceding $J$ intersects $J$, and all preceding right endpoints
lie strictly before the right endpoint of $J$. Hence $W=V\cap J$
is a nonempty proper subinterval of $J$ made up of whole blocks
(Figure~\ref{fig:thin-negative-fibers}(a)).
By minimality of $J$, $\underline b(W)\le y(W)\le\overline b(W)$.

Since $y(V)+y(J)=y(I)+y(W)$, submodularity of $\overline b$ on the
upper fiber and supermodularity of $\underline b$ on the lower fiber give
\begin{equation}\label{eq:5-3}
 0<\Delta_V<\Delta_V+\Delta_J
 \le \Delta_I+\Delta_W\le s .
\end{equation}
If $D\cap V=\varnothing$, then $y(V)$ and $\Delta_V$ are integers.
If $D\subseteq V$, then $y(I\setminus V)$ is an integer and
$y(V)=\Delta-y(I\setminus V)$ gives $\Delta_V\in s+\Z$ on either fiber.
Neither $\Z$ nor $s+\Z$ meets $(0,s)$. Thus $D\cap V$ and $D\setminus V$
are both nonempty, which also rules out $\sigma(y)<0$ when $|D|=1$.

Fix a proper union $V$ of consecutive blocks with $D\cap V$ and
$D\setminus V$ nonempty, and choose $i\in D\cap V$ and $j\in D\setminus V$.
Eliminate $y_j$ using $y(I)=\Delta$ and fix $y_k$ for $k\in D\setminus\{i,j\}$.
By \eqref{eq:source-fiber-coordinates}, the remaining integral is with respect
to $\mathrm{d}y_i$, so the measure of the bold segment in
Figure~\ref{fig:thin-negative-fibers}(b) is the length of its horizontal projection.
As a function of $y_i$, $\Delta_V$ has slope $1$ on the upper fiber and $-1$
on the lower fiber. Thus $0<\Delta_V<s$ restricts $y_i$ to an interval of
length at most $s$. Requiring $0<y_j<c_j$ can only shorten this interval.

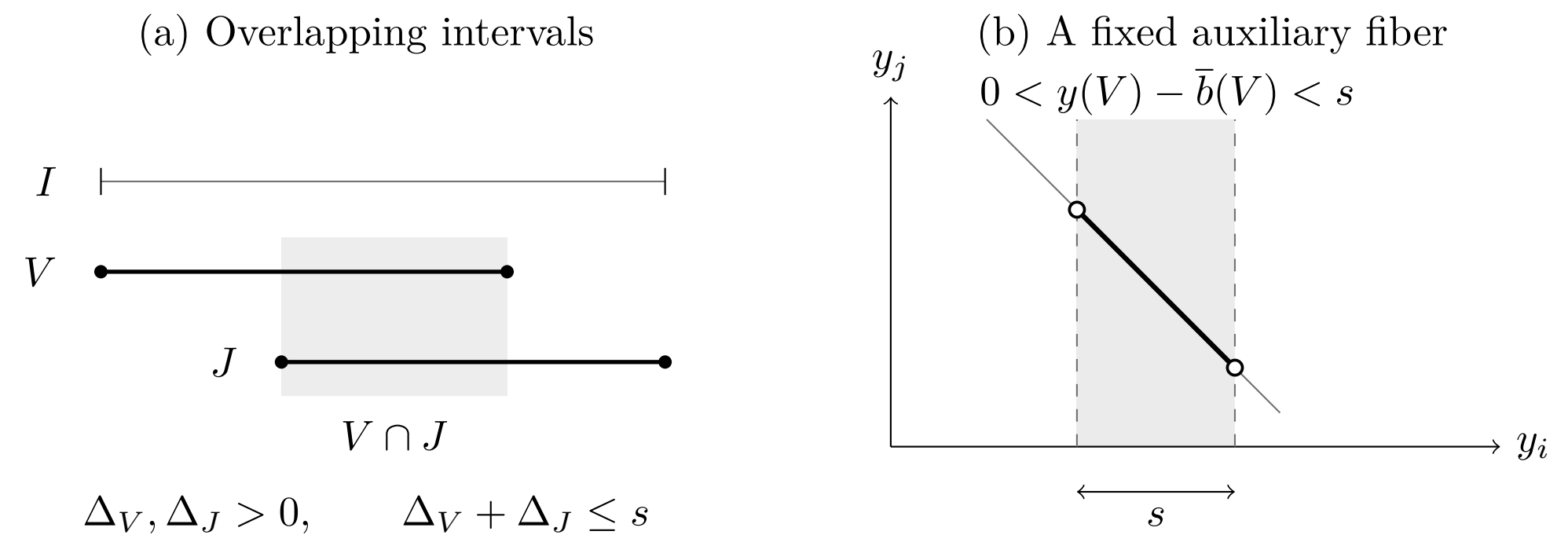
\begin{figure}[htbp]
\centering
\begin{tikzpicture}[x=0.95cm,y=0.95cm,font=\small]
 \node at (3,3.65) {(a) Overlapping intervals};
 \node at (10.5,3.65) {(b) A fixed auxiliary fiber};
 \fill[black!7] (2.25,0.45) rectangle (4.25,1.85);
 \node[anchor=east] at (0.4,2.35) {$I$};
 \draw[black!55] (0.65,2.35)--(5.65,2.35);
 \draw (0.65,2.23)--(0.65,2.47) (5.65,2.23)--(5.65,2.47);
 \node[anchor=east] at (0.4,1.55) {$V$};
 \draw[line width=1pt] (0.65,1.55)--(4.25,1.55);
 \fill (0.65,1.55) circle (1.6pt);
 \fill (4.25,1.55) circle (1.6pt);
 \node[anchor=east] at (2,0.75) {$J$};
 \draw[line width=1pt] (2.25,0.75)--(5.65,0.75);
 \fill (2.25,0.75) circle (1.6pt);
 \fill (5.65,0.75) circle (1.6pt);
 \node at (3.25,0.1) {$V\cap J$};
 \node at (3,-0.6) {$\Delta_V,\Delta_J>0,\qquad \Delta_V+\Delta_J\le s$};
 \fill[black!8] (9.3,0) rectangle (10.7,2.9);
 \draw[->] (7.65,0)--(13.05,0) node[right] {$y_i$};
 \draw[->] (7.65,0)--(7.65,3.1) node[above] {$y_j$};
 \draw[dashed,black!55] (9.3,0)--(9.3,2.9)
                            (10.7,0)--(10.7,2.9);
 \node at (10.1,3.15) {$0<y(V)-\overline b(V)<s$};
 \draw[black!55] (8.5,2.9)--(11.1,0.3);
 \draw[line width=1.2pt] (9.3,2.1)--(10.7,0.7);
 \draw[fill=white,line width=0.7pt] (9.3,2.1) circle (1.8pt);
 \draw[fill=white,line width=0.7pt] (10.7,0.7) circle (1.8pt);
 \draw[<->] (9.3,-0.4)--(10.7,-0.4)
   node[midway,below] {$s$};
\end{tikzpicture}
\caption{The upper fiber $y(I)=h_++s$.}
\label{fig:thin-negative-fibers}
\end{figure}

Integrating over $0<y_k<c_k$ for $k\in D\setminus\{i,j\}$ bounds the
fiber measure of the slab $0<\Delta_V<s$ by
$s\prod_{k\in D\setminus\{i,j\}}c_k$.
For each source, the slabs indexed by the proper block unions $V$ with
$D\cap V$ and $D\setminus V$ nonempty cover
$\{y:\sigma(y)<0\}$. By subadditivity, summing the slab bounds over
the finitely many consecutive block unions $V$ and then over the
finitely many sources gives the stated bound on each fiber.
\end{proof}

For $x\in[0,1]$, apply Lemma~\ref{lem:boundary-limits} to
$\boldsymbol{\psi}_I$ at the cutoffs $h_+$ and $h_-$.
With the auxiliary functions $\lambda_s^\pm$ constructed there,
the coefficientwise limits give
\begin{align}
 \delta_I^+(t,x)&=(1-x)\lim_{s\downarrow0}
   \frac{K_{\boldsymbol{\psi}_I,\lambda_s^+}(0,h_++s)}{t},
   \label{eq:5-1}\\
 \delta_I^-(t,x)&=x\lim_{s\downarrow0}
   \frac{K_{\boldsymbol{\psi}_I,\lambda_s^-}(0,h_--s)}{t}
   \quad\text{if }h_-\ge1. \label{eq:5-2}
\end{align}
When $h_-=0$, the sum defining $\delta_I^-$ in \eqref{eq:2-10} is empty.

\begin{proposition}\label{prop:positive-increments}
Suppose that $\delta_J^\pm(t,x)\in\R_{\ge0}[t]$ for every proper
interval $J\subsetneq I$ with $|J|\ge2$ and every
$x\in[0,1]$. Then, for every $x\in[0,1]$, the boundary increments
$\delta_I^-$ and $\delta_I^+$ computed from $\boldsymbol{\psi}_I$ satisfy
\begin{equation}\label{eq:5-5}
 \delta_I^\pm(t,x)\in\R_{\ge0}[t].
\end{equation}
\end{proposition}

\begin{proof}
Fix $x\in[0,1]$ and an integer $k\ge0$. For $0<s<1$, take
$(\lambda,\Delta)=(\lambda_s^+,h_++s)$ or, when $h_-\ge1$,
$(\lambda,\Delta)=(\lambda_s^-,h_--s)$.
For a source $\mathfrak s\in\mathcal S_I$ with blocks $K_1,\ldots,K_q$,
write $v_J$ for the entry position of block $J$ at initial position $u=0$.
A continuous singleton supplies a factor $t$ in $\upsilon(t,y)$.
After cancelling one such factor, the singleton weights in \eqref{eq:2-13}
have nonnegative coefficients since $0\le\lambda(v_J)\le1$.
Each correction block $J$ is a proper subinterval of $I$, so the
hypothesis gives $\delta_J^\pm(t,\lambda(v_J))\in\R_{\ge0}[t]$.
Thus $\upsilon(t,y)/t\in\R_{\ge0}[t]$.

For fixed $\mathfrak s$, the coefficient $[t^{k+1}]\upsilon(t,y)$ is a
polynomial in $\lambda(v_{K_1}),\ldots,\lambda(v_{K_q})$ and is bounded
on $[0,1]^q$. Since $\mathcal S_I$ is finite and \eqref{eq:4-3} gives
$|\sigma(y)|\le 2^{|I|-1}$, there is a common bound $M_{I,k}$ for
$|\sigma(y)|[t^{k+1}]\upsilon(t,y)$, independent of $x,s,y$.
It follows that
\[
 [t^k]\left(\frac{\sigma(y)\upsilon(t,y)}t\right)
 \ge -M_{I,k}\mathbf1_{\{\sigma(y)<0\}}.
\]
Integrate over $y(I)=\Delta$ for each source and sum over $\mathcal S_I$.
By the kernel formula in Lemma~\ref{lem:exact-inventory}\ref{part:inventory-before}
at $u=0$ and the measure bound in Lemma~\ref{lem:thin-negative-fibers},
\[
 [t^k]\frac{K_{\boldsymbol{\psi}_I,\lambda}(0,\Delta)}{t}
 \ge -M_{I,k}\sum_{\mathfrak s\in\mathcal S_I}
      \mu_{\Delta}\bigl(\{y:\sigma(y)<0\}\bigr)
 \ge -M_{I,k}c_Is.
\]
As $s\downarrow0$, the lower bound tends to zero. The coefficientwise
limits in \eqref{eq:5-1} and \eqref{eq:5-2}, multiplied by $1-x$ and $x$,
respectively, therefore give $[t^k]\delta_I^+(t,x)\ge0$ and, for $h_-\ge1$,
$[t^k]\delta_I^-(t,x)\ge0$.
\end{proof}

\section{Proof of the main theorem}\label{sec:main-proof}

For intervals $J$ of length two, the hypothesis on proper subintervals in
Proposition~\ref{prop:positive-increments} is vacuous. Induction on $|J|$
using that proposition gives $\delta_J^\pm(t,x)\in\R_{\ge0}[t]$ for
every interval $J\subseteq N$ of length at least two and every $x\in[0,1]$.

By Example~\ref{ex:running-singleton}, $\boldsymbol\varphi_{\{i\}}=\boldsymbol\alpha_{c_i}$
is a complete positive state on $[0,c_i]$.
Fix an interval $I$ of length at least two, an auxiliary function $\lambda$
and an initial position $u$. Each source weight is a product of block weights.
For a singleton entered at $v$, the possible weights $\lambda(v)$,
$1-\lambda(v)$ and $t$ belong to $\R_{\ge0}[t]$ since $0\le\lambda(v)\le1$.
For a correction block $J$ entered at $v$, the induction above gives
$\delta_J^\pm(t,\lambda(v))\in\R_{\ge0}[t]$.
The closed guard indicators and the fiber measure $\mu_\Delta$ are nonnegative
and independent of $t$. Hence the sums and fiber integrals in
Lemma~\ref{lem:exact-inventory}\ref{part:inventory-after} give
$\varphi_{I,j}(t,\lambda(u))\in\R_{\ge0}[t]$ for $j\in\Z_{\ge0}$ and
$K_{\boldsymbol\varphi_I,\lambda}(u,u+\Delta)\in\R_{\ge0}[t]$ for positive
noninteger $\Delta$.

Write $h_-=\underline b(I)$ and $h_+=\overline b(I)$.
The two truncations in \eqref{eq:4-1} set $\varphi_{I,j}=0$ for
$j\notin[h_-,h_+]$. For positive noninteger $\Delta\notin[h_-,h_+]$,
\eqref{eq:double-cut-density} gives $K_{\boldsymbol\varphi_I,\lambda}(u,u+\Delta)=0$.
At integer or nonpositive displacements the kernel is zero by
Definition~\ref{def:canonical-kernel}.

To evaluate $\varphi_{I,j}(t,x)$ at $x=1$ or $x=0$, prescribe
$\lambda(u)=x$ as permitted after Definition~\ref{def:auxiliary-function}.
Only sources without continuous branches enter the evaluation weight.
Their jumps are integers, so each block entry $v$ lies in $u+\Z$ and
$\lambda(v)=x$ by periodicity. By \eqref{eq:2-13} and the identities
$\delta_J^+(t,1)=\delta_J^-(t,0)=0$ from \eqref{eq:2-10}, a block $J$
with nonzero weight must have jump $\underline b(J)$ at $x=1$ and
$\overline b(J)$ at $x=0$. Thus for any nonzero evaluation term with
block partition $\pi$, the bounds on its total jump $j$ give
\[
\begin{aligned}
 x=1 &\qquad h_-\le j=\sum_{J\in\pi}\underline b(J)\le h_-,\\
 x=0 &\qquad h_+\ge j=\sum_{J\in\pi}\overline b(J)\ge h_+,
\end{aligned}
\]
where the last inequalities follow from supermodularity of $\underline b$
and submodularity of $\overline b$, respectively. Hence only the component
at $j=h_-$ can be nonzero at $x=1$, and only the component at $j=h_+$
can be nonzero at $x=0$, as required by \eqref{eq:2-12}.
This completes the proof of Proposition~\ref{prop:state-closure}.

\begin{proof}[Proof of Theorem~\ref{thm:main}]
Suppose first that $|E|\ge2$, that $\dim\mathcal P=|E|-1$, and that
$\min_{x\in\mathcal P}x_i=0$ for every $i\in E$. The convolution
$\boldsymbol{\varphi}\coloneqq\boldsymbol{\varphi}_N\boldsymbol{\alpha}_{c_e}$
represents $\beta_N(t,z)\beta_{\{e\}}(t,z)$, whose coefficient polynomial
$\zeta_r$ was identified with $i(\mathcal P,t)$ in Section~\ref{subsec:normal-forms}.
Since $x(N)=r-x_e$ and $x_e$ ranges over $[0,c_e]$ on $\mathcal P$,
\eqref{eq:1-1} gives $\underline b(N)=r-c_e$ and $\overline b(N)=r$.
Proposition~\ref{prop:state-closure} and Example~\ref{ex:running-singleton}
show that $\boldsymbol{\varphi}_N$ and $\boldsymbol{\alpha}_{c_e}$ are complete
positive states on $[r-c_e,r]$ and $[0,c_e]$, respectively. By
Lemma~\ref{lem:positive-products}, $\boldsymbol{\varphi}$ is a complete
positive state on $[r-c_e,r+c_e]$. Since $r-c_e<r$, the endpoint condition
\eqref{eq:2-12} gives $\varphi_r(t,1)=0$. Writing $\delta_r^-$ for the lower
correction of $\boldsymbol{\varphi}$ at $h=r$, we obtain from
\eqref{eq:boundary-coefficient}
\begin{equation}\label{eq:6-1}
 i(\mathcal P,t)=\zeta_r(t)=\delta_r^-(t,1).
\end{equation}
Since $r\ge c_e\ge1$, the lower-limit formula \eqref{eq:boundary-lower-limit}
applies at $h=r$. With $\lambda_s^-$ chosen as in the proof of
Lemma~\ref{lem:boundary-limits} for $x=1$, we obtain
\begin{equation}\label{eq:direct-terminal-limit}
 i(\mathcal P,t)
 =\lim_{s\downarrow0}
   \frac{K_{\boldsymbol{\varphi},\lambda_s^-}(0,r-s)}t
\end{equation}
coefficientwise in $t$. For $0<s<1$, formula \eqref{eq:2-14} gives the
factor $t$ in $K_{\boldsymbol{\varphi},\lambda_s^-}(0,r-s)$, and
Definition~\ref{def:complete-positive-state}\ref{cond:state-positivity}
gives nonnegativity of its coefficients.
Hence, for every integer $k\ge0$,
\[
 [t^k]\,i(\mathcal P,t)
 =\lim_{s\downarrow0}[t^{k+1}]K_{\boldsymbol{\varphi},\lambda_s^-}(0,r-s)\ge0 .
\]

For a general $\mathcal P$, each point factor in \eqref{eq:ehrhart-product}
has Ehrhart polynomial $1$. Every other factor $\mathcal P_a$ has dimension
$|C_a|-1$, so subtracting its integer coordinate minima reduces it to
the case just proved without changing $i(\mathcal P_a,t)$.
\end{proof}

\begin{remark}\label{rem:coordinate-truncation}
We show that, for $0<s<1$ and positive integers $t$ with $st\in\Z$,
the kernel quotient in \eqref{eq:direct-terminal-limit} counts the lattice
points $w\in t\mathcal P$ satisfying $w_e\ge st$.
For the normalized polypositroid $\mathcal P$, write $s=m/q$ in lowest terms and set
$\mathcal P_s=\mathcal P\cap\{x_e\ge s\}$. We claim that $q\mathcal P_s$
is integral and that
\[
 \frac{K_{\boldsymbol{\varphi},\lambda_s^-}(0,r-s)}{t}
 =i(q\mathcal P_s,t/q).
\]
Each new vertex of $q\mathcal P_s=q\mathcal P\cap\{x_e\ge m\}$ lies on an
edge of $q\mathcal P$ crossing $x_e=m$. The primitive edge direction $\mathbf e_i-\mathbf e_j$ has $e$-coordinate
$\pm1$ on a crossing edge. Thus moving from a lattice endpoint to height
$m$ adds an integer multiple of $\mathbf e_i-\mathbf e_j$ and gives a
lattice point.

Take $t$ to be a positive multiple of $q$. The proof of
Lemma~\ref{lem:boundary-limits} at $h=r$ and $x=1$ gives the kernel
quotient as $\sum_{0\le j<r}(\mathsf T_{r-j-s}\varphi_j)[0,1]$.
For $0\le j<r$, put $a=(r-j-s)t$, a positive integer. For $b\ge1$,
\eqref{eq:uniform-wall-count} and \eqref{eq:negative-binomial-coefficient} give
\[
 (\mathsf T_{r-j-s}x^b)[0,1]
 =[z^a](1-z)^{-b}.
\]
For $b=0$, the divided difference of $1$ and $[z^a]1$ both vanish.
By linearity, $(\mathsf T_{r-j-s}\varphi_j)[0,1]$ is the coefficient
of $z^{(r-s)t}$ in $z^{jt}\varphi_j(t,(1-z)^{-1})$. The terms with
$j\ge r$ in \eqref{eq:2-1} have no powers below $jt>(r-s)t$.
Summing over $0\le j<r$ therefore gives the kernel quotient as
$[z^{(r-s)t}]\bigl(\beta_N(t,z)\beta_{\{e\}}(t,z)\bigr)$.

By \eqref{eq:1-3}, each $y\in t\mathcal P_N\cap\Z^N$ lifts uniquely to
$w\in t\mathcal P\cap\Z^E$ with $w_e=rt-y(N)$. To reach exponent
$(r-s)t$, the term $z^{y(N)}$ must pair with $z^{w_e-st}$ in
$\beta_{\{e\}}(t,z)$. Since $st$ is an integer and $0\le w_e\le c_e t$,
the exponent $w_e-st$ belongs to $\{0,\ldots,c_e t\}$ precisely when
$w_e\ge st$. Hence the coefficient equals $\#(t\mathcal P_s\cap\Z^E)$.
At every positive multiple $t=qn$ this count is $i(q\mathcal P_s,n)$,
so the polynomials $K_{\boldsymbol{\varphi},\lambda_s^-}(0,r-s)/t$
and $i(q\mathcal P_s,t/q)$ agree at infinitely many values of $t$.
Consequently, \eqref{eq:direct-terminal-limit} gives
\[
 i(\mathcal P,t)=\lim_{\substack{s=m/q\downarrow0\\ \gcd(m,q)=1}}
 i(q\mathcal P_s,t/q)
\]
coefficientwise in $t$.
\end{remark}

\appendix
\section{Pushforward and product measures}\label{app:measure-conventions}

Let $\tau:D\to H$ be a Borel map between Borel subsets of Euclidean
spaces, and let $\mu$ be a finite Borel measure on $D$. The
\emph{pushforward} of $\mu$ under $\tau$ is the measure $\tau_*\mu$
on $H$ defined by \({(\tau_*\mu)(B)=\mu(\tau^{-1}(B))}\) for every
Borel set $B\subseteq H$.
For every bounded Borel function \({\psi:H\to\R}\),
\begin{equation}\label{eq:pushforward-integral}
 \int_H\psi(x)\,\mathrm{d}(\tau_*\mu)(x)
 =\int_D\psi(\tau(y))\,\mathrm{d}\mu(y).
\end{equation}

\begin{definition}\label{def:occupation-measure}
Let $a<v$ be real numbers and let $\tau:(a,v)\to H$ be a Borel
measurable map into a Borel subset $H$ of a Euclidean space.
The \emph{occupation measure} of $\tau$ on $(a,v)$ is the finite
Borel measure $\rho$ on $H$ defined by
\[
 \rho(B)=\int_a^v\mathbf1_B\bigl(\tau(s)\bigr)\,\mathrm{d}s
\]
for every Borel set $B\subseteq H$, where $\mathbf1_B$ denotes the
indicator function of $B$.
\end{definition}

For the auxiliary function $\lambda$ of Section~\ref{subsec:integral-operators},
let $\rho_{a,v}$ denote its occupation measure on $(a,v)$.
With $\tau=\lambda$ and Lebesgue measure on $(a,v)$,
\eqref{eq:pushforward-integral} gives
\begin{equation}\label{eq:occupation-values}
 \int_{[0,1]}\psi(x)\,\mathrm{d}\rho_{a,v}(x)
 =\int_a^v\psi(\lambda(s))\,\mathrm{d}s
\end{equation}
for every bounded Borel function $\psi:[0,1]\to\R$.

Let $\rho$ be a finite Borel measure on $H$. For an integer $j\ge1$,
write $H^j=H\times\cdots\times H$, with $j$ factors.
The \emph{product measure} $\rho^{\otimes j}$ on $H^j$ is determined by
\[
 \rho^{\otimes j}(B_1\times\cdots\times B_j)
 =\prod_{i=1}^j\rho(B_i)
\]
for Borel sets $B_1,\ldots,B_j\subseteq H$.
For $j=0$, $\rho^{\otimes0}$ assigns mass one to the singleton $H^0$.

\section*{Acknowledgments}
This work was supported by the National Natural Science Foundation of China
(Nos.~12325112, 12288101).

\Needspace{6\baselineskip}

\bigskip
\Needspace{8\baselineskip}
{\small
\noindent
\textsc{Y. Li, School of Mathematical Sciences, Peking University,
Beijing 100871, P.R. China.}
\par\nobreak
\noindent
Email address: \texttt{\href{mailto:liyaao@stu.pku.edu.cn}{liyaao@stu.pku.edu.cn}}

\medskip
\noindent
\textsc{S. Shao, CAPT, LMAM and School of Mathematical Sciences,
Peking University, Beijing 100871, P.R. China.}
\par\nobreak
\noindent
Email address: \texttt{\href{mailto:sihong@math.pku.edu.cn}{sihong@math.pku.edu.cn}}
\par}
\end{document}